\documentclass[11pt]{article}
\usepackage{amsfonts}
\usepackage{graphicx}
\usepackage{amscd}
\usepackage{amsmath}
\newlength{\dinwidth}
\newlength{\dinmargin}
\def\be{\begin{equation}}
\def\ee{\end{equation}}
\def\ben{\begin{displaymath}}
\def\een{\end{displaymath}}
\def\baa{\begin{eqnarray}}
\def\eaa{\end{eqnarray}}

\def\ba{\begin{array}}
\def\ea{\end{array}}
\makeatletter \@addtoreset{equation}{section} \makeatother

\newtheorem{theorem}{Theorem}

\newtheorem{remark}{Remark}
\newtheorem{lemma}{Lemma}

\def\dvadva{{{{\cal M}}(\Z_2)}}
\def\dom{\Omega(\Z_2)}

\def\Z{{\bf Z}}
\def\L{\cal L}
\def\pmo{{\bf e}}

\def\p{\partial}
\def\be{\begin{equation}}
\def\ee{\end{equation}}
\def\la{\label}
\def\th{\Theta}

\def\L{{\cal L}}

\def\p{\partial}
\def\C{{\mathbb C}}
\def\B{{\bf B}}
\def\F{{\cal F}}
\def\M{{\cal M}}
\def\D{{\cal D}}

\def\S{{\cal H}}

\def\G{{\cal G}}
\def\mod{{\rm Sp}(4,\Z)}
\def\U{{\frak{U}}}
\def\a{{\bf x}}
\def\b{{\bf y}}

\def\g{\gamma}

\def\s{\sigma}

\def\th{\vartheta}
\def\Th{\Theta}

\def\e{\epsilon}
\def\DD{{\rm det}\Delta}
\def\fo{{f}} 
\def\f{\frac}

\begin{document}
\title{Extremal properties of the determinant of the 
Laplacian in the Bergman metric  on the moduli space of
genus two Riemann surfaces}

\author{C.Klein$^1$, A. Kokotov$^2$, D. Korotkin$^2$}

\maketitle

\begin{center}
$^1$Max-Planck Institute for Mathematics in the Sciences, \\
        Inselstra\ss e 22,       04103 Leipzig, Germany\\
$^2$Department of Mathematics and Statistics, Concordia University,\\
7141 Sherbrooke West, Montreal H4B 1R6, Quebec,  Canada
\end{center}
\vskip0.5cm
{\bf Abstract.} We study extremal properties 
of the determinant of the Laplacian
in the Bergman metric on the moduli space of compact genus two Riemann surfaces. 
By a combination of analytical and numerical methods we identify four non-degenerate 
critical points of this function and compute the signature of the 
Hessian at these points.
The curve with the maximal number of automorphisms (the Burnside curve) turns out to be the 
point of the absolute maximum. Our results agree with the mass formula 
for virtual Euler characteristics 
of the moduli space. A similar analysis is performed for
Bolza's strata of symmetric Riemann surfaces of genus two. 
\vskip0.5cm
\tableofcontents

\section{Introduction}

The study of extremal properties of various functionals  related to
 Riemann surfaces  attracted the attention of
many researchers during the last 25 years (see \cite{Sarnak} and other
papers of the same volume  for an introduction to the subject). 
The functionals studied so far are related to both 
purely geometrical aspects of a Riemann surface 
as the function $syst$ (the length of the
shortest closed geodesics, see \cite{Schmutz} and references therein), and 
to spectral aspects as the minimal eigenvalue of the
corresponding Laplace operator \cite{Jacobson}, or the (appropriately
regularized) determinant of the Laplacian \cite{OPS}. The spectral
characteristics are determined not only by the conformal structure on a
Riemann surface, but also by the choice of a metric within a given
conformal class. As was shown in \cite{OPS} (for a short proof
see \cite{KKRicci}),  the determinant of the Laplacian $\det \Delta$
within a given conformal class of compact Riemann surfaces
of fixed genus $g\geq 0$ and fixed volume takes its
maximum for the  metric of constant curvature;
therefore this determinant 
 defines a natural functional on the moduli space of
Riemann surfaces. 


So far the study of
extremal properties of ${\rm det}\Delta$ in metrics of constant
curvature did not go beyond the genus one case, where it is possible
to prove that the tori with periods equal to $i$ and $e^{{\pi i}/{3}}$
are critical points of $\DD$ on the moduli space: 
$e^{{\pi i}/{3}}$ is the maximum, and $i$ is the saddle point. 
The proof of the vanishing of the gradient of $\DD$ at these two points 
is based on the existence of subgroups of the modular group
leaving these points invariant \cite{OPS} (i.e., the 
corresponding elliptic surfaces have non-trivial automorphisms
groups); thus the gradient of any
modular invariant functional, not only $\DD$, vanishes at these two points.

In genus two and higher the analysis of extremal properties of any
functional becomes decisively more complicated. Nonetheless some
explicit results are known in genera two and three (it was recently
proved in \cite{Jacobson} that the minimal eigenvalue of the Laplacian in
genus two is maximal on the Burnside curve equipped with a singular
metric with cone singularities); the same holds for the function $syst$  \cite{Schmutz}
(the Jacobian of  the Burnside
curve also defines a lattice in $\C^2$ corresponding to the densest
sphere packing \cite{Sarnak}). In \cite{Schmutz} $syst$ was 
used as a topological Morse function to define a cell decomposition of the moduli spaces;
in the simplest cases (in genera $0,1,2$ with a small number of punctures) the
analysis of the critical points of $syst$ allowed to reproduce known values of
virtual Euler characteristics of the moduli spaces \cite{HZ} via the mass formula.

Unfortunately any reasonably complete
treatment of extremal properties of the determinant of the Laplacian in the
Poincar\'e metric in genus two seems to be out of reach at the
moment due to the absence of efficient numerical algorithms for the
computation of this determinant (perhaps such an algorithm can be 
developed by extending results of \cite{PoRo}
to the full moduli space). 

The goal of this paper is to study extremal properties of another
smooth functional on the moduli space of genus two Riemann surfaces with good boundary behavior - the
function
\be
\F= ({\rm det}\Im \B)^{5/2}\prod_{s=1}^{10}|\Th[\beta_s](\B)|
\la{F}
\ee
where $\B$ is the matrix of $b$-periods of a compact Riemann surface $\L$ computed with respect to
some canonical basis of cycles $(a_j,\,b_j)$ 
(expression (\ref{F}) does not depend on the choice of this basis);
the product is taken over all ten even theta characteristics $\beta_s$;  $\Th[\beta_s](\B)$
is the theta constant corresponding to the characteristic $\beta_s$ (in this paper we denote by
$\Th$ the genus two theta function; the genus one theta functions are denoted by $\th$).

The function $\F$ is related to the determinant of the Laplace operator as follows:
$$
{\rm det}\Delta_B= C\; \F^{1/3}\;,
$$
where $C$ is a moduli-independent constant; $\Delta_B$ is the Laplacian in the 
Bergman metric - the metric of volume $1$ given by
\be
g_B:=\sum_{j,k=1}^2 (\Im\B)^{-1}_{jk}v_j \bar{v}_k\;,
\la{Bergmet}
\ee
where the $v_j$ are holomorphic 1-forms on $\L$ normalized by $\oint_{a_j} v_k=\delta_{jk}$. 

The Bergman metric is induced on $\L$ by  the standard flat invariant K\"ahler metric on the Jacobian 
when $\L$ is canonically embedded into its Jacobian via the Abel map.
Since in genus two the theta divisor 
can be biholomorphically mapped to the Riemann surface itself, 
$\F^{-1/3}$ turns out to coincide also with
the analytic torsion of the theta divisor equipped with the metric induced 
by the same K\"ahler metric on the Jacobian \cite{Yoshikawa}.

Being considered as a function on the upper Siegel half-space $\S$, the function $\F$
coincides with the Petersson norm $||\Delta_2||=({\rm det}\Im \B)^{5/2}|\Delta_2|$
of the Siegel cusp form $\Delta_2:=\prod_{s=1}^{10}\Th[\beta_s]$. Finally the 
function $\F$ essentially coincides with the
genus two Mumford measure \cite{Knizhnik}.

On the boundary of the moduli space $\M$ of genus two Riemann surfaces $\F$ vanishes;
therefore  $-\log\F$ is the proper function on $\M$; moreover,  all of its critical points are 
non-degenerate (which is equivalent to the non-degeneracy of critical points of $\F$ itself), 
which together with the boundary behavior $-\log\F\to +\infty$ would allow 
the study of topological properties of $\M$. 

In this paper we study critical points 
of $\F$ on $\M$ by using a combination of analytical and numerical tools. 
We prove that any smooth function on the upper Siegel half-space $\S$ 
invariant with respect to $\mod$ has critical points corresponding to 
three curves with
large automorphism groups: the Burnside curve $y^2=x(x^4-1)$, the $D_6$ curve $y^2=x^6-1$ and
the $\Z_5$ curve $y^2=x^5-1$ (this result is in particular valid for ${\rm det} \Delta$ in
the Poincar\'e metric!). Further numerical analysis of $\F$ in  Gottschling's
fundamental domain  $\S/\mod$  shows that all of these critical points 
are non-degenerate, that the Burnside curve gives the absolute 
maximum of $\F$, and that the signature of the Hessian at these three points 
is equal to $(0,6)$, $(3,3)$ and $(2,4)$ respectively.
In addition the numerical analysis shows the existence of a fourth critical 
point -- a curve from the $D_3$ family
 $y^2=(z^3-1)(z^3-r^3)$ with
$$
r=0.22373907612077\ldots
$$
This point is also non-degenerate, and the signature of the Hessian 
there is $(1,5)$.

Though the presence of the fourth critical point appears somewhat 
unexpected at first, its existence is
predicted by the mass formula \cite{Schmutz} for orbifold Euler 
characteristics of $\M$. The mass formula states that  the Euler characteristic
is equal to the sum over all critical Riemann surfaces $\L_i$ of $(-1)^{c_i}/\{\#{Aut(\L_i)}\}$
where $c_i$ is the index of the critical point;  $\#{Aut(\L_i)}$ is the order of 
the group of automorphisms of $\L_i$. One can easily check that the mass formula
immediately implies the existence of the fourth critical curve.

We perform a similar analysis on three smaller moduli spaces of genus two curves with fixed groups
of automorphisms: the space of (complex) dimension two of curves with 
a $\Z_2\times\Z_2$
symmetry group, and one-dimensional spaces of curves with $\Z_2\times D_2$ and $\Z_2\times D_3$
symmetry groups; it turns out that $\F$ does not have any new critical points on these
sub-spaces. 
In particular, we show that the moduli space of  $\Z_2\times D_2$ 
curves in coordinates given by the
one remaining parameter in the matrix of $b$-periods can be identified with 
the factor of the upper half-plane by the modular
group $\Gamma_0(2)+$; similarly, the  moduli space of  $\Z_2\times D_3$ curves can be identified with the factor of the
upper half-plane by the  modular
group $\Gamma_0(3)+$. We discuss the meaning of the mass formula for these  subspaces
of symmetric curves.

The paper is organized as follows: In Section 2 we recall what is known about extremal 
properties of ${\rm det} \Delta$ on the moduli space of genus one Riemann surfaces.
In Section 3 we summarize necessary facts on the description of the full moduli space,
as well as of its symmetric strata, in terms of matrices of $b$-periods. Although most of the
facts presented here are well-known, we did not find some of them (about the link 
between $D_2$ and $D_3$
moduli spaces with subgroups $\Gamma_0(2)+$ and $\Gamma_0(3)+$ of the modular group) 
in the existing literature.  In Section 4 we prove that the three curves with large automorphism groups
are stationary points of $\F$, as well as of any other smooth function on
$\S$ invariant under the action of the $\mod$ group. Numerical 
analysis shows the existence of an additional
critical point on $\M$ and the non-degeneracy of all four critical points.
The Burnside curve turns out to be the global maximum of $\F$ on $\M$; 
 we compute the signature of the Hessian for the other critical points.
In Section 5 a similar analysis is performed for the 
three strata of $\M$ with given degree of symmetry.
In Section 6 we outline the relationship of our results to the 
computation of virtual Euler characteristics of $\M$ and  of its symmetric subspaces.

In the sequel we shall use the following  notation for the root of unity: $\e_k=e^{2\pi i/k}$;
the modular group $SL(2,\Z)$ will be denoted by $\Gamma$.

\section{Summary of the genus one case}

In genus one the Bergman metric (\ref{Bergmet}) coincides with the metric of constant curvature - the
 flat metric of volume $1$ given by
${dz d\bar{z}}/{\Im\sigma}$
on the torus with periods $1$ and $\sigma$. The determinant of the Laplacian 
(acting on functions, i.e., sections of the trivial line bundle) in this metric is, up to a multiplicative constant,  given by
the expression \cite{RaySinger}:
\be
\fo(\sigma)= C (\Im \sigma)^{1/2}|\eta(\sigma)|^2\;,
\la{Fgen1}
\ee
where $\eta=[\th_1^{\prime}]^{1/3}$ is the Dedekind eta-function; $C$ is a constant independent
 of the moduli. Due to the Jacobi formula $\th_1^{\prime}=i\th_2\th_3\th_4$ ($\th_j$, $j=2,3,4$ are the standard theta-constants) the expression (\ref{Fgen1})
is a straightforward analog of the genus two expression (\ref{F}).

The function $\fo$ (\ref{Fgen1}) is modular invariant, real and positive.
Moreover it vanishes on the boundary of the moduli space, when the 
torus degenerates and when
$\sigma$ tends to $+i\infty$ (or any other point related to $+i\infty$ by a modular transformation).
The function $\fo$ has the following obvious symmetry with respect to 
reflections at the imaginary axis:
\be
\fo(-\bar{\s})=\fo(\s)
\la{mirgen1}
\ee

The extremal properties of the function $\fo$  (\ref{Fgen1}) are well-known (see \cite{OPS,Sarnak}):
\begin{theorem}
The points $\sigma=i $ and $\sigma = e^{\pi i/3}$ (as well as all 
points obtained from these two points by modular transformations) are stationary 
points of the function (\ref{Fgen1}); the point
$e^{\pi i/3}$ is the point of the absolute maximum of $\fo$.
\end{theorem}

 Both   points $i$ and $e^{\pi i/3}$  are orbifold points 
(with cone angle $\pi$ and $2\pi/3$ respectively) of the moduli space which 
can be obtained by an appropriate identification of the boundary of
the fundamental domain $\Omega$ of the group $\Gamma$:
there are non-trivial subgroups of the modular group leaving these two points invariant.
Introduce the  standard generators of  $\Gamma$:
\be
t=\left(\ba{cc} 1 & 1 \\
                0 & 1 \ea\right)\;,\hskip1.0cm
s=\left(\ba{cc} 0 & -1 \\
                1 & 0 \ea\right)\;.
\la{genG}
\ee
The stationary subgroup of order 2 leaving the point $\sigma=i$ invariant is generated by the element  $s$.
The stationary subgroup of order 3 leaving the point $\sigma= e^{\pi i/3}$ invariant is 
generated by the product $st$.  

\begin{figure}[hbt]
\centering 
   \includegraphics[width=10cm]{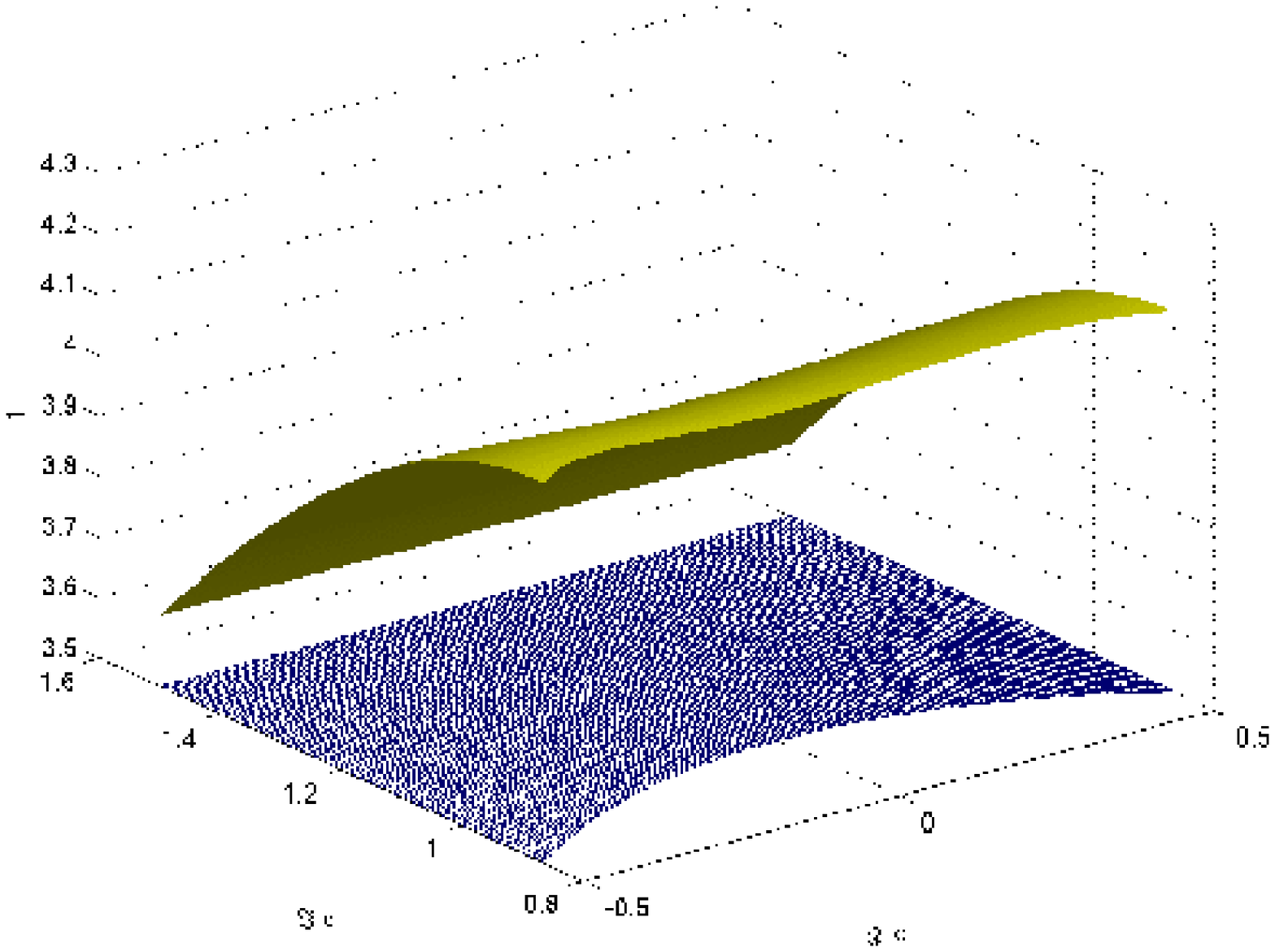}
   \caption{Plot of $\fo(\s)$ in the fundamental domain $\Omega$}
   \label{figfOmega}
\end{figure}

The proof of the stationarity \cite{OPS} 
of $\fo(\s)$ at these two points works equally well for {\it any}
smooth function on the upper half-plane invariant under the modular group;
below we prove a similar statement in the genus two case. 
It was proved analytically in \cite{OPS} that the point $e^{\pi i/3}$ is the 
local maximum of $f(\s)$; it can be shown numerically that actually the point $e^{\pi i/3}$ is the absolute maximum, and 
$i$ is the saddle point of $\fo$,  (see Fig.~\ref{figfOmega} for the plot
of $\fo$ in the fundamental domain $\Omega$ of the modular group).

\begin{figure}[hbt]
\centering 
  \includegraphics[width=10cm]{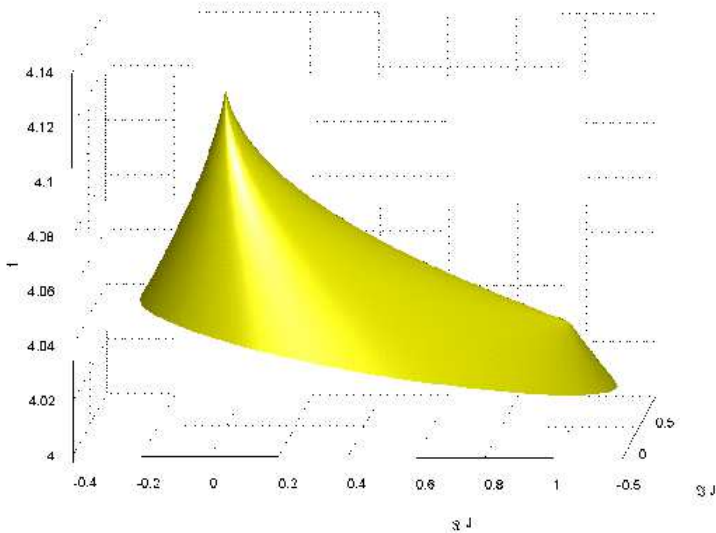}
    \caption{Plot of $\fo$ as a function of the $J$-invariant}
   \label{figfJ}
\end{figure}

It is  instructive to plot $\fo$ also as a function of the $J$-invariant:
\be
J(\sigma):=\f{(\th_2^8+\th_3^8+\th_4^8)^3}{54\th_2^8\th_3^8\th_4^8}\;,
\ee
which maps the fundamental domain $\Omega$ onto the whole complex plane.
The maximum of $\fo$ in Fig.~\ref{figfJ} is achieved at $J(e^{\pi i/3})=0$ (the highest peak);
the spike at $J(i)=1$ is the saddle point. 
The spikes of the function $\fo(J)$ at these points appear since the change of
coordinates $\sigma\to J$ is degenerate there.

\section{Moduli space of genus two Riemann surfaces}

\subsection{The Siegel fundamental domain: Gottschling's description}

The moduli space of Riemann surfaces $\M$ is covered (with  branching) by the Torelli space 
${\cal T}$ which is the space of marked Riemann surfaces, i.e., the space of
pairs (the Riemann surface $\L$ and the canonical basis of cycles on $\L$). The Torelli space 
is still not simply-connected; its fundamental group is called the Torelli 
group; the  universal 
covering of the Torelli space is the Teichm\"uller space. 
It is important for us that the covering of the Torelli space by the 
Teichm\"uller space is unramified, i.e., that any
analytic function on the Torelli space remains analytic whilst 
being lifted to the Teichm\"uller space.

 Consider the upper Siegel half-space 
$\S=\{z\in M(2, {\mathbb C})\ :\ z=z^t; \ \Im z\geq 0\}\;.$
Assigning to any marked Riemann surface its matrix of \( b \)-periods 
computed in a given basis, we
get the map from the Torelli space to $\S$; the change of a canonical basis 
of cycles on a 
given Riemann surface corresponds to an $\mod$ transformation of $\S$. 
The image of this map does not coincide with the whole space
$\S$ since the matrix of $b$-periods
of a non-degenerate Riemann surface of genus two can never be 
diagonal or equivalent to diagonal
up to a modular transformation. Denote by $\D\subset \S$ the set of matrices 
which are either 
diagonal or can be transformed to diagonal form by a modular transformation. 
In genus two, when the independent entries of the matrix of $b$-periods 
can be used as
local coordinates on the moduli space, the space  $\S\setminus \D$
can be identified with the Torelli space. The factor of  $\S\setminus \D$ 
by the action of 
$\mod$ can be identified with the moduli space $\M$ of genus two Riemann surfaces.

We can also first factorize  $\S$ by the action of $\mod$.  This gives the 
Siegel-Gottschling fundamental domain $\G:=\S/\mod$; taking out from $\G$ points 
lying in $\D$,  we get  the space $\G\setminus \D$
which also coincides with the moduli space $\M$.

The fundamental domain $\G$ can be described by 25 inequalities on
the matrix entries of $\B$ \cite{Gottschling}. 
To describe these conditions we introduce real and imaginary 
parts of the independent components of $\B$:
$\mathbf{B}_{11}=x_{1}+iy_{1}$, $\mathbf{B}_{12}=x_{2}+iy_{2}$ and $\mathbf{B}_{22}=x_{3}+iy_{3}$. 

Then $\G$ is defined by the following set of inequalities:
\begin{itemize}
\item 
Conditions restricting the range of $\{x_i\}$ and $\{y_i\}$:
\begin{equation}    
|x_{i}|\leq \frac{1}{2}, \quad i=1,2,3,\quad     
 y_{i}\geq\frac{1}{2}\sqrt{3}, \quad i=1,3, \quad y_{2}\geq0    
\label{eq:num1}.
\end{equation}

\item the Minkowski ordering condition: 
\begin{equation}    y_{1}\geq 2y_{2}, \quad y_{3}\geq y_{1}    
\label{mink},
\end{equation}
\item The following set of  19 inequalities:
\begin{equation}    
|\mathbf{B}_{11}|\geq 1,\quad |\mathbf{B}_{22}|\geq 1,\quad    |\mathbf{B}_{11}+\mathbf{B}_{22}-2\mathbf{B}_{12}+\pmo|\geq1,    \label{gott1}\end{equation}and 
\begin{equation}    |\det (\mathbf{B}+S)|\geq 1    
\label{gott2},
\end{equation}
where $S$ are the matrices
\begin{equation}    
\begin{split}    
\begin{pmatrix}        0 & 0  \\        0 & 0    
\end{pmatrix},\quad     
\begin{pmatrix}     \pmo & 0  \\            0 & 0       
\end{pmatrix},\quad     \begin{pmatrix}     0 & 0  \\       0 & \pmo    
\end{pmatrix},\quad     \begin{pmatrix}     \pmo & 0  \\            0 & \pmo    
\end{pmatrix},\quad     \\    \begin{pmatrix}       \pmo & 0  \\            0 & -\pmo   
\end{pmatrix},\quad     \begin{pmatrix}     0 & \pmo  \\            \pmo & 0    
\end{pmatrix},\quad     \begin{pmatrix}     \pmo & \pmo  \\         \pmo & 0    
\end{pmatrix},\quad     \begin{pmatrix}     0 & \pmo  \\            \pmo & \pmo 
\end{pmatrix}    \label{gott3}    
\end{split}
\end{equation}
and $\pmo=\pm 1$.
\end{itemize}

\subsection{Curves with non-trivial automorphisms: Bolza's classification}

Any  genus two Riemann surface is biholomorphically equivalent to  an algebraic curve
defined by an equation:
$
y^2=f(z)
$
where $f(z)$ is a polynomial of degree 5 or 6.
The hyperelliptic involution on this curve maps any point  $(y,z)$  
to $(-y,z)$; this involution generates the hyperelliptic 
symmetry group $\Z_2$. A generic genus two curve does not have
other automorphisms. 
If a curve has a larger automorphism group $Aut$, the hyperelliptic $\Z_2$ subgroup  always 
turns out to be a normal subgroup
of $Aut$. All symmetric curves can be stratified according to the type of 
the reduced automorphism group
$Aut/\Z_2$ which can be one of the following 6 types \cite{Bolza}:
$\Z_2$, $D_2$, $D_3$, $S_4$, $D_6$ or $\Z_5$
(where $D_i$ are dihedral groups).

\subsubsection{ Group $\Z_2$}

 Any Riemann surface $\L$ from this family can be represented by the equation
\be 
y^2=(z^2-1)(z^2-r_1^2)(z^2-r_2^2)
\la{Z2} \ee 
where $r_1,\,r_2\in\C$.
In addition to the  hyperelliptic involution there is an involution $\mu$ acting as 
$\mu:\; z\to -z$ on the curve (\ref{Z2}).
\begin{figure}[htb]
    \centering 
     \includegraphics[width=10cm]{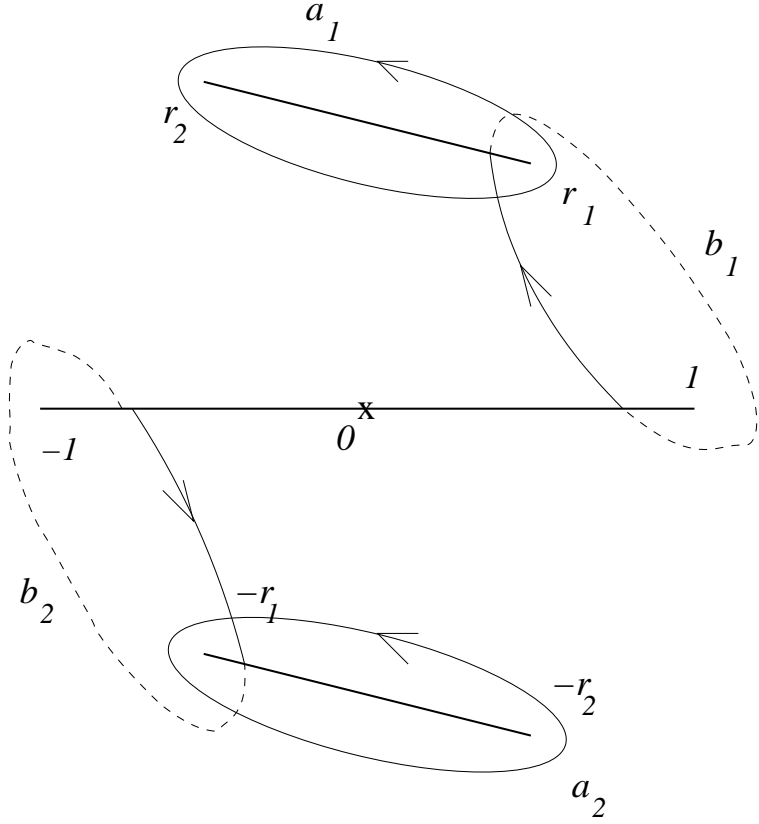}
    \caption{Canonical basis of cycles on a curve with $\Z_2\times\Z_2$ symmetry}
    \label{Z2fig}
\end{figure}
The basic cycles on $\L$ chosen as shown in Fig.~\ref{Z2fig} transform as
follows under the action of the involution $\mu$: 
\be
a_1^{\mu}=a_2\;\hskip0.7cm b_1^\mu=b_2\;,
\la{abmu}
\ee 
or equivalently,
 \be
\left(\ba{cccc} b_1^\mu\\ b_2^\mu\\ a_1^\mu\\a_2^\mu\ea\right)=
T^\mu \left(\ba{cccc} b_1\\ b_2\\ a_1\\a_2 \ea\right) 
\ee 
where the  $\mod$ matrix $T^\mu$ is given by 
\be 
T^\mu= 
\left(\ba{cccc}           0 & 1 & 0 & 0 \\
                          1 & 0 & 0 & 0 \\
                          0 & 0 & 0 & 1 \\
                          0 & 0 & 1 & 0 \ea \right)
\ee
Writing the matrix $T^\mu$ in block form:
\be
T^\mu=\left(\ba{cc} A^\mu & B^\mu\\
                    C^\mu & D^\mu\ea\right)
\ee
where $ B^\mu= C^\mu=0$;  $ A^\mu= D^\mu= \left(\ba{cc} 0 & 1\\
                                                        1 & 0\ea\right)$,
and taking into account that, on one hand, the matrix of $b$-periods transforms 
under the action of any symplectic transformation as follows:
\be
\B\to \B^\mu := (A^\mu\B+ B^\mu)(C^\mu\B+ D^\mu)^{-1}\;,
\ee
and that, on the other hand, the matrix of $b$-periods must remain invariant 
under the action of the biholomorphic transformation mapping the 
canonical basis of cycles to  the new one, we conclude that $\B=\B^\mu$. 
This is equivalent to
$$
\B=\left(\ba{cc} 0 & 1\\
                 1 & 0\ea\right) \B 
\left(\ba{cc} 0 & 1\\
              1 & 0\ea\right)\;;
$$
thus the matrix $\B$ can be parametrized as follows:
\be
\B=\f{1}{2}\left(\ba{cc}\a+\b & \a-\b\\
                        \a-\b & \a+\b\ea\right)
\la{z2prym}
\ee
If we factorize $\L$ with respect to the involution $\mu$, we get an elliptic curve $\L_0$ with period 
$\b$; $\a$ is the period of elliptic Prym variety corresponding to the (ramified) covering 
$\L\to\L_0$.

The description of the moduli space $\dvadva$ was given in \cite{Schiller} in terms
of the variables $\tau_1=\b$ and $\tau_2=-1/\a$. 
Below we describe this construction in terms of $\a$ and $\b$ themselves
which makes it slightly more transparent.
 Denote by 
$\Gamma(2)$ the main congruence subgroup  of $\Gamma$ consisting 
of matrices $\gamma$ such that $\gamma\equiv I \; ({\rm mod} \;2)$ ($I$ is the unit matrix). 

The natural idea of \cite{Schiller} is to reduce the action of the full 
group $\mod$ on matrices of the form (\ref{z2prym}) to a correlated action of $\Gamma$
on $\a$ and $\b$.

\begin{lemma}
Assume that the basis of canonical cycles $(b_{1},b_{2},a_1,a_2)$ satisfying 
(\ref{abmu})
is transformed by a matrix $T\in\mod$ to another basis  $(b'_1,b'_2,a'_1,a'_2)$ 
satisfying the same relation:
\be
a_1'^\mu= a_2'\;,\hskip0.8cm b_1'^\mu= b_2'\;.
\la{abmupr}\ee
Then the action of the matrix $T$ on the matrix of $b$-periods (\ref{z2prym}) gives rise to an action of
two elements $\g_1,\g_2\in\Gamma$, 
such that $\g_1\g_2^{-1}\in\Gamma(2)$, on $\a$ and $\b$, respectively. The matrix $T$ is expressed in terms of $\g_{1,2}$ as follows: if
\be
\g_i=\left(\ba{cc} k_i & l_i \\
                   m_i & n_i \ea\right)\;,\hskip1.0cm i=1,2
\la{gi}
\ee
then
\be
T=\f{1}{2}\left(\ba{cccc} k_1+k_2 & k_1-k_2  & l_1+l_2 &  l_1-l_2 \\
                          k_1-k_2 & k_1+k_2  & l_1-l_2 &  l_1+l_2 \\
                          m_1+m_2 & m_1-m_2  & n_1+n_2 &  n_1-n_2 \\
                          m_1-m_2 & m_1+m_2  & n_1-n_2 &  n_1+n_2 \ea\right)
\la{Tz2z2}
\ee
\la{ghZ2}
\end{lemma}
{\it Proof.} The additional symmetries of the matrix  $T\in\mod$ (\ref{Tz2z2}) 
follow from the assumptions (\ref{abmu}),
(\ref{abmupr}) about the behavior of the basic cycles under the 
action of the involution $\mu$. The 
equivalence of the action of the matrix $T$  on the matrix (\ref{z2prym}) 
to the action of the two elements
$\g_{1,2}$ (\ref{gi}) separately on $\a$ and $\b$ follows by direct computation.
Finally the condition that all matrix entries of $T$ (\ref{Tz2z2}) are integer 
is equivalent to the condition that all matrix entries of the matrix $\g_1-\g_2$ 
are even. A simple computation
using the conditions ${\rm det}\g_i=1$ shows that this is equivalent 
to the condition 
$\g_1\g_2^{-1}\in\Gamma(2)$  (which is more natural from the point of view of the group structure); 
it is straightforward to check that this condition
is invariant under the natural group operation in $\Gamma\times\Gamma$.

$\square$

As was proved in \cite{Schiller}, to get the full subgroup of 
$\mod$ preserving the form (\ref{z2prym}) of the matrix of $b$-periods, 
one has to add to the subgroup (\ref{Tz2z2}) 
one more transformation, given by
\be
T= 
\left(\ba{cccc}           1 & 0 & 0 & 0 \\
                          0 & -1 & 0 & 0 \\
                          0 & 0 & 1 & 0 \\
                          0 & 0 & 0 & -1 \ea \right)
\la{addit}
\ee
which maps the canonical basis $(b_1,b_2,a_1,a_2)$ to the basis $b_1'=b_1\;, b_2'=-b_2$,  
and $a_1'=a_1\;, a_2'=-a_2$; obviously, this transformation destroys the symmetry (\ref{abmu}) 
of the canonical basis of cycles under the involution $\mu$. The 
action of the transformation (\ref{addit}) on the
matrix of $b$-periods (\ref{z2prym}) is very simple: it changes the sign of 
the off-diagonal terms,
i.e., it corresponds to the interchange of $\a$ and $\b$: $\a'=\b$, $\b'=\a$.
Let us denote by $S_2$ the two-element permutation group generated by this 
transformation.

The theorem proved in \cite{Schiller} can now be reformulated as follows:
\begin{theorem}\la{funz2}
The moduli space $\dvadva$ can be represented as the following factor:
\be
\dvadva= {\cal S}(\Z_2)/G\;,
\la{dvadva}
\ee
where
\be
{\cal S}(\Z_2)=H\times H\setminus \{(\a,\b)\;|\; \a=\g\b\;,\a,\b\in H\;,\; \g\in \Gamma(2)\}
\la{fd22}
\ee
($H$ is the upper half-plane), and where the group $G$ is defined as follows:
\be
G= G_0  S_2
\la{grG}
\ee
where $G_0  S_2 $  is the semi-direct product of two groups  and 
$G_0$ is the following normal subgroup
of $G $:
\be
G_0=\{(\g_1,\g_2)\in \Gamma\times\Gamma\; |\; \g_1\g_2^{-1}\in\Gamma(2)\}
\la{G_0}
\ee
\end{theorem}
\vskip0.3cm
The subspace  $\{(\a,\b)\;|\; \a=\g\b\;,\a,\b\in H\;,\; \g\in \Gamma(2)\}$ is taken out of $H\times H$
since it consists of matrices which are $\mod$ equivalent to diagonal ones. The appearance of the
group $G$ is explained above; Lemma \ref{ghZ2} shows how to construct the natural group homomorphism 
$f$ from $G$ to $\mod$.

The space ${\cal S}(\Z_2)$ (\ref{fd22}) is called in \cite{Schiller} the ``special Torelli space'', and the group
$G$ (\ref{grG}) the ``special Torelli group''.

The non-trivial part of the proof of this theorem (for which we refer the reader to \cite{Schiller})
 is to show that
no $\mod$ transformation  from the complement of the image of the homomorphism $f$
preserves the matrix of $b$-periods  (\ref{z2prym}).

\begin{figure}[htb]
    \centering 
    \includegraphics[width=10cm]{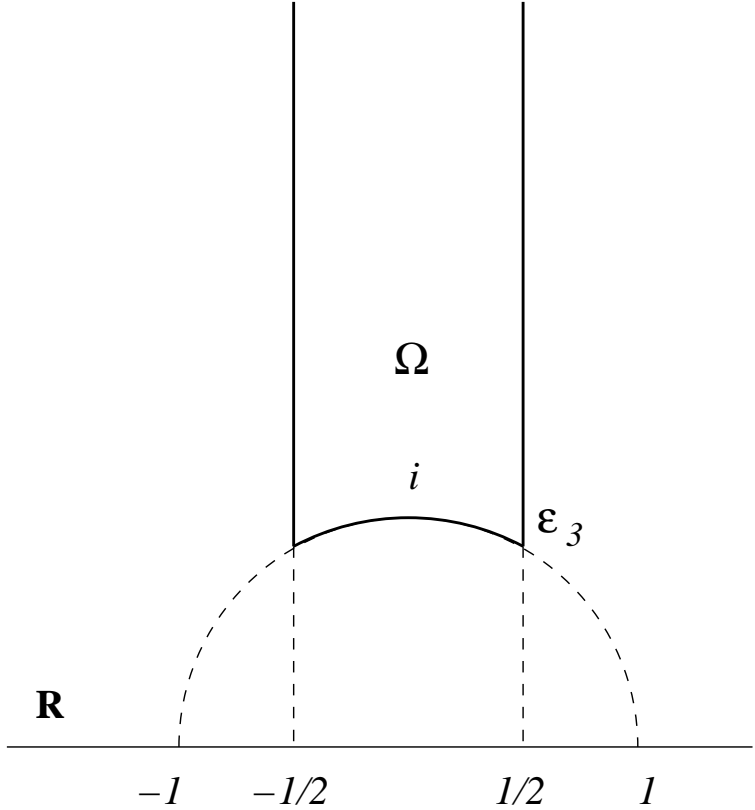}
    \caption{Fundamental domain $\Omega$ of the modular group $\Gamma$ }
    \label{Omegafig}
\end{figure}

\begin{figure}[htb]
    \centering 
     \includegraphics[width=10cm]{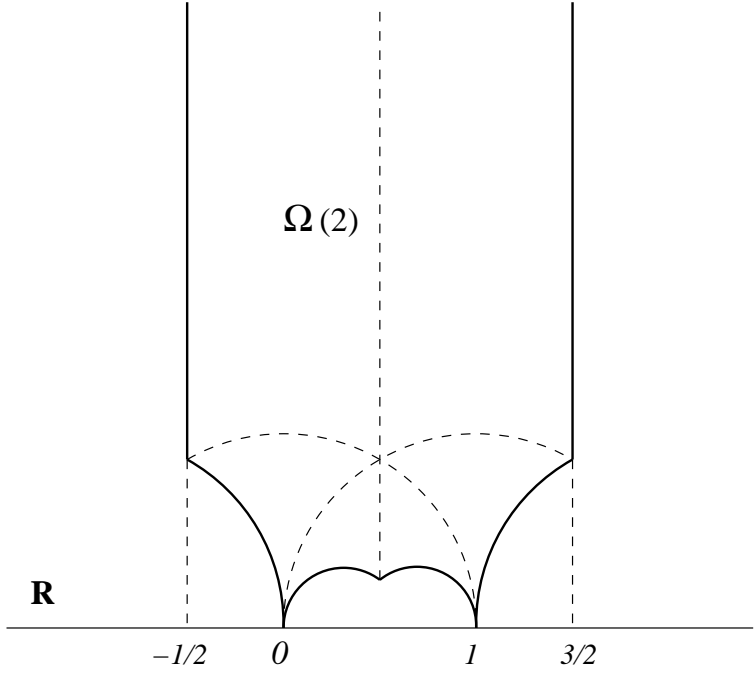}
    \caption{Fundamental domain $\Omega(2)$ of the subgroup $\Gamma(2)$}
    \label{Omega2fig}
\end{figure}

The structure of the fundamental domain $\dvadva$ (\ref{dvadva})
 is rather non-trivial due to the necessity to 
take into account the subgroup $S_2$ interchanging $\a$ and $\b$. Therefore for our
 subsequent numerical analysis we introduce the ``bigger'' factor space
$\dom=(H\times H)/G_0$, which can be easily 
described.
\begin{lemma}\la{simfunz2}
The fundamental domain for the action of the group $G_0$ on $H\times
H$ can be chosen to be $\Omega\times\Omega(2)$, where $\Omega$ is the
standard fundamental domain of the group $\Gamma$ (Fig.\ref{Omegafig}), 
and $\Omega(2)$ (Fig.\ref{Omega2fig})  is
the fundamental domain of the subgroup $\Gamma(2)$, consisting of six
copies of $\Omega$.
\end{lemma}
{\it Proof} is simple. To prove that any point $(\a,\b)\in H\times H$
can be mapped to the inside of  $\Omega\times\Omega(2)$ by some transformation
from $G_0$, we first identify a $\g_1$ such that $\g_1\a\in \Omega$. As the
second step we find a transformation $\g\in \Gamma(2)$ such that
$\g(\g_1\b)\in \Omega(2)$. Obviously the transformation  mapping
$(\a,\b)$ to $(\g_1\a,\,\g\g_1\b)$ belongs to $G_0$ since $\g_1
(\g\g_1)^{-1}= \g^{-1} \in \Gamma(2)$.

Suppose now that some transformation $(\g_1,\,\g_2)\in G_0$ maps some
point $(\a,\b)\in \Omega\times\Omega(2)$ to another point $(\a',\b')\in
\Omega\times\Omega(2)$.
Since $\g_1\in\Gamma$, it must be the identity element (since $\Omega$ is a
fundamental domain for $\Gamma$). Therefore $\g_2\in\Gamma(2)$. Since $\g_2$
maps a point $\b\in \Omega(2)$ to another point $\b'\in\Omega(2)$,
$\g_2$ must also be the identity transformation.

$\square$

\subsubsection{Group $D_2$}

Curves of this family, which forms a subfamily of the two-parametric family (\ref{Z2}),
can be represented by the equation
\be
y^2=z(z^2-1)(z^2-r^2)
\la{D2}
\ee
where $r\in \C$. In addition to the hyperelliptic involution, 
there are two more involutions: $\mu_1:\, z\to  -z$ and $\mu_2:\,z\to r/z$
on the curve  (\ref{D2});
the order of the full symmetry group $D_2\times\Z_2$ equals $8$.
\begin{figure}[htb]
    \centering 
     \includegraphics[width=10cm]{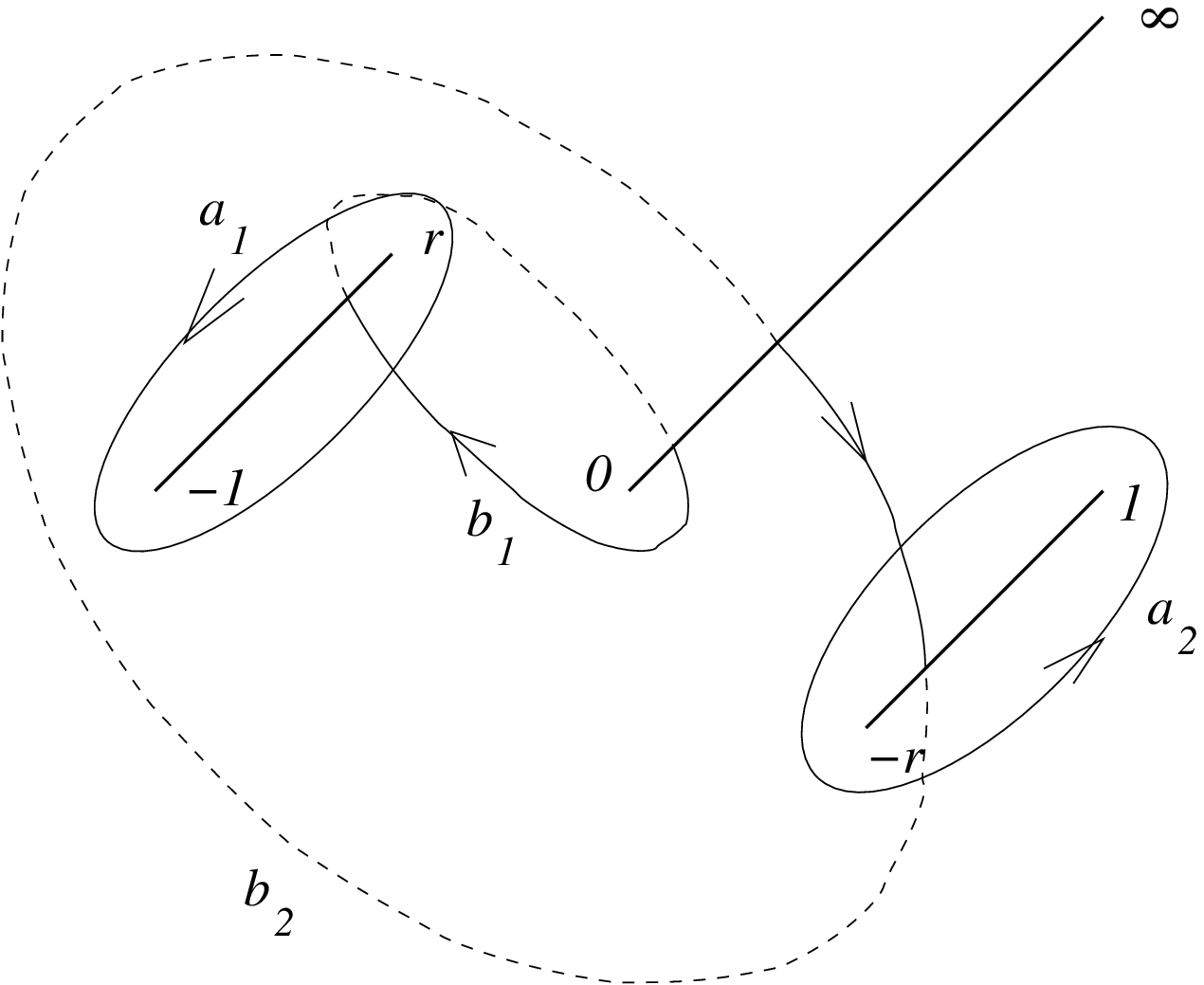}
    \caption{Canonical basis of cycles on a curve with $D_2\times \Z_2$ symmetry}
  \label{D2fig}
\end{figure}
Choose the branch cuts and basic cycles $(a_i, b_i)$ on $\L$ as shown 
in Fig.~\ref{D2fig}.  The involutions $\mu_{1,2}$ act
on this basis as follows:
\be
b_1^{\mu_1}=b_2-a_1\;, \hskip0.5cm b_2^{\mu_1}=-b_1+a_2\;, \hskip0.5cm
a_1^{\mu_1}=a_2\;, \hskip0.5cm a_2^{\mu_1}=-a_1\;,
\ee
and
\be
b_1^{\mu_2}=b_2\;, \hskip0.5cm b_2^{\mu_2}=b_1\;,\hskip0.5cm
a_1^{\mu_2}=a_2\;, \hskip0.5cm a_2^{\mu_2}=a_1\;;
\ee
 therefore the action of $\mu_1$ and $\mu_2$ on
the vector $(b_1,\,b_2,\,a_1,\,a_2)^t$ is given by the following $\mod$ matrices:
\be
T^{\mu_1}=\left(\ba{cccc} 0 & 1 & -1 & 0 \\
                          -1 & 0 & 0 & 1 \\
                          0 & 0 & 0 & 1 \\
                          0 & 0 & -1 & 0 \ea \right)\;;\hskip0.7cm
T^{\mu_2}=
\left(\ba{cccc}           0 & 1 & 0 & 0 \\
                          1 & 0 & 0 & 0 \\
                          0 & 0 & 0 & 1 \\
                          0 & 0 & 1 & 0 \ea \right)\;,
\ee
respectively.
Invariance of the matrix of $b$-periods under the $\mod$ transformation defined by $T^{\mu_1}$ 
implies the following structure of the matrix  in the chosen basis:
\be
\B=\left(\ba{cc} \s & 1/2 \\
                 1/2 & \s \ea\right)
\la{BD2}
\ee
where $\s\in \C$; $\Im \s >0$. Invariance of the matrix of $b$-periods 
under the action of $T^{\mu_2}$ does not impose any other restriction.

It will be convenient for us to work in terms of the parameter
$\a=\s+1/2$, i.e., to represent the matrix $\B$ in the form
\be
\B=\left(\ba{cc} \a-1/2 & 1/2 \\
                 1/2 & \a-1/2 \ea\right)
\la{BD22}
\ee
for $\Im\a>0$.
The moduli space $\M(D_2)$ of $D_2$ curves can be realized as the upper complex
half-plane in the variable $\a=\s+1/2$ factorized with respect to the action of 
the group $\mod$ on matrices of the form (\ref{BD22}).

The analogue of Schiller's theorem \ref{funz2} for
the $D_2$  family looks as follows.
Let us introduce the modular group $\Gamma_0(2)$ which is the
subgroup of $SL(2,\Z)$ consisting of
matrices whose $(21)$ entry is even.
The minimal subgroup of $GL(2,{\bf Q})$ containing $\Gamma_0(2)$ and
the element
$$
 \left(\ba{cc} 0 & -1 \\
               2 & 0 \ea\right)
$$
is called $\Gamma_0(2)+$.
 The group $\Gamma_0(2)+$ is generated by the following two elements \cite{Ford,Chris}:
\be
\g_1= \left(\ba{cc} 1 & 1 \\
                    0 & 1 \ea\right)\;\hskip0.7cm {\rm and} \hskip0.7cm
\g_2= \left(\ba{cc} 2 & -1 \\
                    2 & 0 \ea\right)\;.
\la{genD2}
\ee

\begin{theorem}\la{lemg2}
The moduli space $\M(D_2)$ can be represented as the following factor:
\be
\M(D_2)= {\cal S}(D_2)/\Gamma_0(2)+
\ee
where
\be
{\cal S}(D_2)= H\setminus \left\{\gamma \left(\f{1}{2}+\f{i}{2}\right)\;,\;\;\; \gamma\in \Gamma_0(2)+\right\}
\la{SD2}
\ee
\end{theorem}
{\it Proof.} 
The proof can be obtained as a corollary of Schiller's theorem \ref{funz2}. To restrict $\M(\Z_2)$ to
$\M(D_2)$ one should put $\b=\a-1$. It is straightforward to verify that the  subgroup of the 
group $G$ (\ref{grG}) preserving the constraint $\b=\a-1$ coincides with $\Gamma_0(2)+$.

In particular,  the
action of the generators $\g_1$ and $\g_2$ on $\a$ is equivalent to the
action of the following  matrices $T_{1,2}\in \mod$ on the matrix  
(\ref{BD22}), respectively:
\be
f(\g_1):=T_1=\left(\ba{cccc} 1 & 0 & 1 & 0 \\
                              0 & 1 & 0 & 1 \\
                              0 & 0 & 1 & 0 \\
                              0 & 0 & 0 & 1  \ea\right)\;,\hskip1.0cm
f(\g_2):=T_2=\left(\ba{cccc} 0 & -1 & 1 & 0 \\
                             0 & -1 & 0 & 0 \\
                             -1 & -1 & 0 & 0 \\
                             1 & -1 & 1 & -1  \ea\right)\;.
\la{T12D2}
\ee
This defines the group homomorphism $f:\Gamma_0(2)+\to \mod$ on the
generators of $\Gamma_0(2)+$ which extends  to the whole group 
according to the group structure.

It remains to understand which matrices of the form (\ref{BD22}) are $\mod$-equivalent
to diagonal ones. Since it is not obvious to get this information from Schiller's 
theorem \ref{funz2}, we choose an indirect way. 
First, applying the
transformation
$$
\left(\ba{cccc} 0 & 1 & 0 & 0 \\
                0 & 1 & -1 & 0 \\
               -1 & 0 & 0 & 1 \\
               1 & 0 & 0 & 0  \ea\right)
$$
to the matrix (\ref{BD2}) with $\a=1/2+i/2$, we get the diagonal matrix
${\rm diag}(i,\;i)$. Therefore, this vertex of the fundamental domain 
$\Omega_0(2)+$ 
(as well as any point equivalent to this point under
$\Gamma_0(2)+$ transformations), corresponds to a degeneration of a genus two Riemann surface
to the union of two tori.
On the other hand, the moduli space $\M(D_2)$ can be alternatively parametrized by
  $\alpha^2$, where $\alpha$ is
the coefficient of the sextic
$xy(x^4+\alpha x^2 y^2 +y^4)$ (this form of equation can be obtained by a
simple transformation from (\ref{D2})); all values of $\alpha^2$ correspond to 
different points of $\M(D_2)$. In terms of $\alpha^2$ the space 
$\M(D_2)$ looks like the Riemann sphere
with two deleted points: $\alpha^2= \infty$, and $\alpha^2=4$, where the Riemann surface degenerates
\cite{Getzler}.
Since we already know two degeneration points in terms of the parameter $\a$
($\a=i\infty$ and $\a=1/2+i/2$, which correspond to $\alpha^2=\infty$ and $\alpha^2=4$, respectively),
 we can conclude that
these are the only boundary points  of $\M(D_2)$ in $\Omega_0(2)+$.

Therefore, all matrices of the form (\ref{BD2}) which are 
$\mod$-equivalent to diagonal ones,
can be obtained 
from the value $\a=1/2+i/2$ by a $\Gamma_0(2)+$ transformation. 

$\square$ 

\begin{figure}[htb]
    \centering 
     \includegraphics[width=10cm]{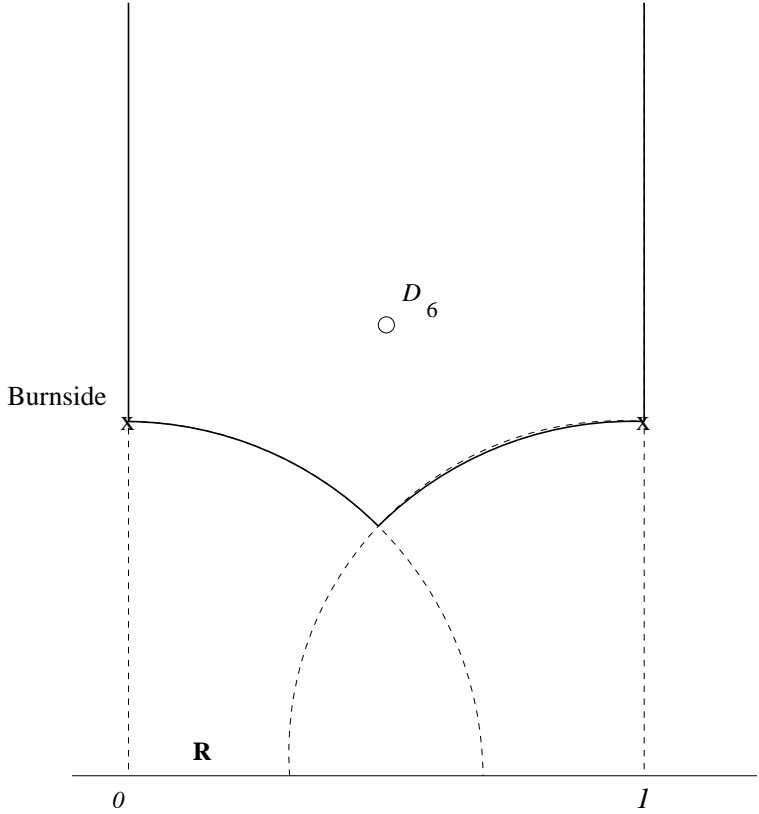}
    \caption{Fundamental domain $\Omega_0(2)+$ with critical  points 
    of $\F$. The bounding 
circles have radius $1/\sqrt{2}$ and centers at $0$ and $1$. The Burnside curve 
corresponds to two values of $\s$ marked by  ${\rm x}$; the point 
marked by a small circle corresponds to the $D_6$ curve}
    \label{Omega200fig}
\end{figure}

Now we can identify the space $\M(D_2)$ with the fundamental domain
$\Omega_0(2)+$ of the group $\Gamma_0(2)+$ (see Fig.\ref{Omega200fig}). The vertical 
lines are identified by the transformation $\g_1$ from (\ref{genD2}); the arcs are identified by 
the transformation $\g_2$ from  (\ref{genD2}).


The space ${\cal S}(D_2)$ (\ref{SD2}) is the natural analog of the Torelli space for the $D_2$ family 
(the ``special $D_2$ Torelli space''); the group $\Gamma_0(2)+$ can be naturally called the ``special
$D_2$ Torelli group'', in analogy to Schiller's terminology for the case of $\M(\Z_2)$.

\subsubsection{Group $D_3$}

This is another subfamily of the two-parametric family (\ref{Z2}) \cite{Bolza}.
Curves admitting this symmetry group also form a one-parametric sub-family of (\ref{Z2}); they
can be represented by the equation
\be
y^2=(z^3-1)(z^3-r^3)
\la{D3}
\ee
where $r\in \C$. In addition to the hyperelliptic involution two more 
generators of the symmetry group are acting on the curve  (\ref{D3}): the element of order three 
$\mu_1: z\to \e_3 z$ and the involution $\mu_2: z\to r/z$;
the order of the full symmetry group $D_3\times\Z_2$ equals $12$.
\begin{figure}[htb]
    \centering 
    \includegraphics[width=10cm]{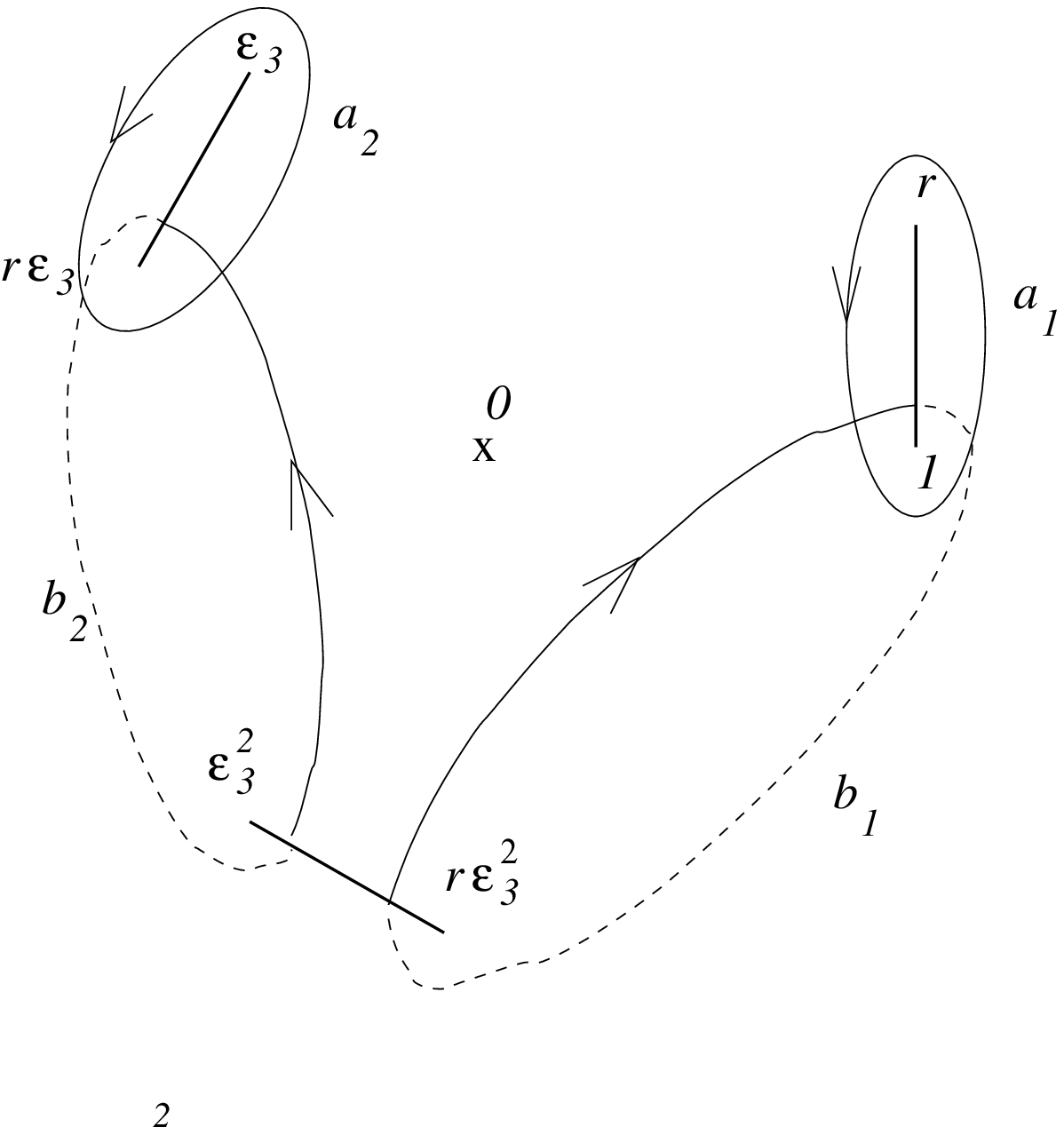}
    \caption{Canonical basis of cycles on a curve with $D_2\times\Z_2$ symmetry}
    \label{D3fig}
\end{figure}
Choose the branch cuts and basic cycles $(a_i, b_i)$ on $\L$ as shown 
in Fig.~\ref{D3fig}.  The symmetries $\mu_{1,2}$ act
on this basis as follows:
\be
b_1^{\mu_1}=b_2-b_1\;, \hskip0.5cm b_2^{\mu_1}=-b_1\;, \hskip0.5cm
a_1^{\mu_1}=a_2\;, \hskip0.5cm a_2^{\mu_1}=-a_1-a_2\;; 
\ee
and 
\be
b_1^{\mu_2}=b_1-b_2\;, \hskip0.5cm b_2^{\mu_2}=-b_2\;,\hskip0.5cm
a_1^{\mu_2}=a_1\;, \hskip0.5cm a_2^{\mu_2}=-a_1-a_2\;;
\ee
therefore the action of $\mu_1$ and $\mu_2$ on 
the vector $(b_1,\,b_2,\,a_1,\,a_2)^t$ is given by the following $\mod$ matrices:
\be
T^{\mu_1}=\left(\ba{cccc} -1 & 1 & 0 & 0 \\
                          -1 & 0 & 0 & 0 \\
                          0 & 0 & 0 & 1 \\
                          0 & 0 & -1 & -1 \ea \right)\;;\hskip0.7cm
T^{\mu_2}=
\left(\ba{cccc}           1 & -1 & 0 & 0 \\
                          0 & -1 & 0 & 0 \\
                          0 & 0 & 1 & 0 \\
                          0 & 0 & -1 & -1 \ea \right)\;,
\ee
respectively.

Invariance of the matrix of $b$-periods under the transformation 
defined by $T^{\mu_1}$ implies its following structure in the chosen basis:
\be
\B=\left(\ba{cc} 2\s & \s \\
                 \s & 2\s \ea\right)
\la{BD3}
\ee
where $\s\in\C$; $\Im\s>0$.  Invariance of the matrix of $b$-periods matrix under the action 
of $T^{\mu_2}$ does not impose any other restriction.

Introduce the subgroup $\Gamma_0(3)$ of $SL(2,\Z)$, which consists of
matrices whose $(21)$ matrix entry $\equiv 0 \;(mod\;\;3)$. The
minimal subgroup of $GL(2,{\bf Q})$ containing both $\Gamma_0(3)$ and the
matrix 
$$
\left(\ba{cc} 0 & -1 \\
              3 & 0 \ea\right)
$$ 
is called $\Gamma_0(3)+$.
The group $\Gamma_0(3)+$ is generated by the following two elements:
\be
\g_1= \left(\ba{cc} 1 & 1 \\
                    0 & 1 \ea\right)\;,\hskip1.0cm
\g_2= \left(\ba{cc} 3 & -1 \\
                    3 & 0 \ea\right)\;;
\la{genG30}
\ee
\begin{theorem}\la{lemg3}
The moduli space $\M(D_3)$ can be represented as the following factor:
\be
\M(D_3)= {\cal S}(D_3)/\Gamma_0(3)+
\ee
where
\be
{\cal S}(D_3)= H\setminus \left\{\gamma 
\left(\f{1}{2}+\f{i}{2\sqrt{3}}\right)\;,\;\;\; \gamma\in 
\Gamma_0(3)+\right\}
\la{SD3}
\ee
\end{theorem}
{\it Proof.} 
The proof can be obtained as a corollary of Schiller's theorem \ref{funz2}. To restrict $\M(\Z_2)$ to
$\M(D_3)$ one should put $\a=3\s$, $\b=\s$.
It is straightforward to verify that the  subgroup of the 
group $G$ (\ref{grG}) preserving the constraint $\a=3\b$ coincides with $\Gamma_0(3)+$ acting on 
$\s(\equiv\b)$.

In particular, the
action of the generators $\g_1$ and $\g_2$  on $\s$ is equivalent to the
action of the following $\mod$  matrices $T_{1,2}\in \mod$ on the matrix of
$b$-periods (\ref{BD3}), respectively:
$$
f(\g_1):=T_1=\left(\ba{cccc} 1 & 0 & 2 & 1 \\
                               0 & 1 & 1 & 2 \\
                               0 & 0 & 1 & 0 \\
                               0 & 0 & 0 & 1  \ea\right)\;,\hskip1.0cm
f(\g_2):=T_2=\left(\ba{cccc}   -2 & 1 & 1 & 0 \\
                               -1 & 2 & 0 & -1 \\
                               -1 & 0 & 0 & 0 \\
                               0 & 1 & 0 & 0  \ea\right)\;,
$$

This defines  the group homomorphism $f:\Gamma_0(3)+\to \mod$ on the
generators of $\Gamma_0(3)+$ which extends  to the whole group 
according to the group structure.

It remains to understand which matrices of the form (\ref{BD3}) are $\mod$- equivalent
to diagonal ones. Again, it is not obvious to get this information from Schiller's 
theorem \ref{funz2}, and we shall use well-known facts about the structure of
$\M(D_3)$ in terms of the coefficients of the sextic.

First, applying the
transformation
$$
\left(\ba{cccc} 0 & 1 & 0 & 0 \\
                0 & 1 & -1 & 0 \\
                -1 & 0 & 0 & 1 \\
                1 & 0 & 0 & 0  \ea\right)
$$
to the matrix (\ref{BD3}) with $\a=1/2+i/2\sqrt{3}$, we get the diagonal matrix
${\rm diag}(-\f{1}{2}+\f{i\sqrt{3}}{2},\;\f{1}{2}+\f{i\sqrt{3}}{2})$. Therefore, this point (as well as any point equivalent to this point under
$\Gamma_0(3)+$ transformations), corresponds to a degeneration of a genus two Riemann surface
to the union of two tori.

On the other hand, the moduli space $\M(D_3)$ can be alternatively parametrized by
  $\alpha^2$, where $\alpha$ is
the coefficient of the sextic
 $x^6+\alpha x^3 y^3-y^6$ (this form of the equation can be obtained by a
simple transformation from (\ref{D3})); all values of $\alpha^2$ correspond to 
different points of $\M(D_2)$. In terms of $\alpha^2$ the space 
$\M(D_3)$ looks like the Riemann sphere
with two deleted points: $\alpha^2= \infty$, and $\alpha^2=-4$, where the Riemann surface degenerates
\cite{Getzler}.
Since we already know two degeneration points in terms of the parameter $\a$
($\a=i\infty$ and $\a=1/2+i/2\sqrt{3}$, which correspond to $\alpha^2=\infty$ and $\alpha^2=-4$, respectively),
 we can conclude that
these are the only  boundary points of $\M(D_3)$ in  $\Omega_0(3)+$.

Therefore, all matrices of the form (\ref{BD3}) which are $\mod$-
equivalent to a  diagonal one,
can be obtained 
from the value $\a=1/2+i/2\sqrt{3}$ by a $\Gamma_0(3)+$ transformation. 

$\square$ 

\begin{figure}[htb]
    \centering 
     \includegraphics[width=10cm]{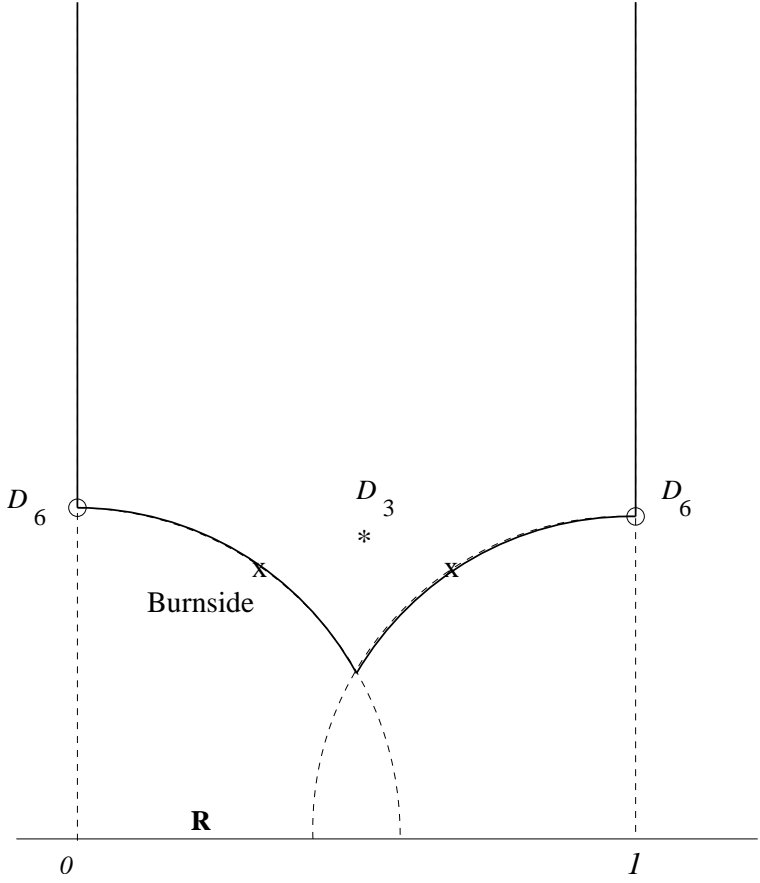}
\caption{Fundamental domain $\Omega_0(3)+$. The bounding 
circles have radius $1/\sqrt{3}$ and centers at $0$ and $1$. The small circles mark the $D_6$ curve,
${\rm x}$ denotes the Burnside curve, and $*$ denotes the critical $D_3$ curve}
    \label{Omega300fig}
\end{figure}

Now we can identify the space $\M(D_3)$ with the fundamental domain
$\Omega_0(3)+$ of the group $\Gamma_0(3)+$ (see Fig.~\ref{Omega300fig}). The vertical 
lines are identified by the transformation $\g_1$ from (\ref{genG30}); the arcs are identified by 
the transformation $\g_2$ from  (\ref{genG30}).


The space ${\cal S}(D_3)$ (\ref{SD3}) is the natural analog of the Torelli space for the $D_3$ family 
(the ``special $D_3$ Torelli space''); the group $\Gamma_0(3)+$ can be naturally called the ``special
$D_3$ Torelli group''.

\subsubsection{ Group $S_4$} The curve admitting this reduced 
symmetry group is defined by the equation
\be
y^2=z(z^4-1)\;;
\la{Bolza}
\ee
it is called the Burnside curve. In 
addition to the hyperelliptic involution the 
generators of the symmetry group of the Burnside curve are  given by $\mu_1:\; z\to i z$\;,\;
$\mu_2:\; z\to (z+1)/(z-1)$ and $\mu_3:\;z\to -\f{1}{z}$. The order of the
full symmetry group $S_4\times\Z_2$ of the Burnside curve equals 48. The Burnside 
curve belongs to both the $D_2$ 
and the $D_3$ 
families.

Choose the branch cuts and basic cycles $(a_i, b_i)$ on $\L$ as shown 
in Fig.~\ref{D2fig} for $r=i$ .  The generators $\mu_{1,2,3}$ act
on this basis as follows:
\be
b_1^{\mu_1}=-b_2+a_1-a_2\;, \hskip0.5cm b_2^{\mu_1}=b_2+a_2\;, \hskip0.5cm
a_1^{\mu_1}=-b_1+b_2-a_1+a_2\;, \hskip0.5cm a_2^{\mu_1}=-b_1-b_2\;,
\ee
\be
b_1^{\mu_2}=-b_1\;, \hskip0.5cm b_2^{\mu_2}=b_1+b_2\;,\hskip0.5cm
a_1^{\mu_2}=b_2-a_1+a_2\;, \hskip0.5cm a_2^{\mu_2}=-b_1+a_2\;;
\ee
and
\be
b_1^{\mu_3}=b_2\;, \hskip0.5cm b_2^{\mu_3}=b_1\;,\hskip0.5cm
a_1^{\mu_3}=a_2\;, \hskip0.5cm a_2^{\mu_3}=a_1\;;
\ee
therefore the action of $\mu_1$, $\mu_2$ and $\mu_3$ on
the vector $(b_1,\,b_2,\,a_1,\,a_2)^t$ is given by the following matrices:
\be
T^{\mu_1}=\left(\ba{cccc} 0 & -1 & 1 & -1 \\
                          0 & 1 & 0 & 1 \\
                          -1 & 1 & -1 & 1 \\
                          -1 & -1 & 0 & 0 \ea \right)\;,\hskip0.7cm
T^{\mu_2}=
\left(\ba{cccc}           -1 & 0 & 0 & 0 \\
                          1 & 1 & 0 & 0 \\
                          0 & 1 & -1 & 1 \\
                          -1 & 0 & 0 & 1 \ea \right)\;,\hskip0.7cm
T^{\mu_3}=
\left(\ba{cccc}           0 & 1 & 0 & 0 \\
                          1 & 0 & 0 & 0 \\
                          0 & 0 & 0 & 1 \\
                          0 & 0 & 1 & 0 \ea \right)\;,
\la{stabB1}
\ee
respectively.
Invariance of the matrix of $b$-periods under the transformations defined 
by $T^{\mu_1}$ and $T^{\mu_3}$ implies the following structure of the matrix 
in the chosen basis:
\be
\B=\left(\ba{cc} -\f{1}{2}+\f{i}{\sqrt{2}} & \f{1}{2} \\
                               \f{1}{2}    & -\f{1}{2}+\f{i}{\sqrt{2}} \ea\right)\;;
\la{Burnside}
\ee
this matrix is also invariant under the transformation defined by  $T^{\mu_2}$.

\subsubsection{Group $D_6$} The curve admitting this reduced symmetry 
group is defined by the equation
\be
y^2=z^6-1\;.
\la{D6}
\ee
In addition to the hyperelliptic involution, there are two more independent
generators of the symmetry group acting on the curve (\ref{D6}): the element
of order six $\mu_1: z\to \e_6 z$ and the involution $\mu_2:\;z\to -1/z$. 
The order of the full symmetry group $D_6\times\Z_2$
equals 24. This curve belongs to both the $D_2$ 
and $D_3$  families.

Choose the branch cuts and basic cycles $(a_i, b_i)$ on $\L$ as shown 
in Fig.~\ref{D3fig} for $r=\e_6$.  The involutions $\mu_{1,2}$ act
on this basis as follows:
\be
b_1^{\mu_1}=-a_1\;, \hskip0.5cm b_2^{\mu_1}=-a_1-a_2\;, \hskip0.5cm
a_1^{\mu_1}=b_1-b_2\;, \hskip0.5cm a_2^{\mu_1}=b_2\;,
\ee
and
\be
b_1^{\mu_2}=b_2\;, \hskip0.5cm b_2^{\mu_2}=b_1\;,\hskip0.5cm
a_1^{\mu_2}=a_2\;, \hskip0.5cm a_2^{\mu_2}=a_1\;;
\ee
therefore 
the vector $(b_1,\,b_2,\,a_1,\,a_2)^t$ transforms under the action of $\mu_1$ and $\mu_2$ 
by the following $\mod$ matrices:
\be
T^{\mu_1}=\left(\ba{cccc} 0 & 0 & -1 & 0 \\
                          0 & 0 & -1 & -1 \\
                          1 & -1 & 0 & 0 \\
                          0 & 1  & 0 & 0 \ea \right)\;;\hskip0.7cm
T^{\mu_2}=
\left(\ba{cccc}           0 & 1 & 0 & 0 \\
                          1 & 0 & 0 & 0 \\
                          0 & 0 & 0 & 1 \\
                          0 & 0 & 1 & 0 \ea \right)\;,
\la{stabB2}
\ee
respectively.
Invariance of the matrix of $b$-periods under 
transformations defined by $T^{\mu_1}$ and $T^{\mu_2}$ 
implies the following structure of the matrix 
in the chosen basis:
\be
\B=\f{i}{\sqrt{3}}\left(\ba{cc} 2 & 1 \\
                                1 & 2 \ea\right)\;.
\la{BD6}
\ee

\subsubsection{ Group $\Z_5$}

The curve admitting this reduced symmetry group is defined by the equation
\be
y^2=z^5-1\;.
\la{Z5}
\ee
The generator of the $\Z_5$ symmetry is given by $\mu: \;  z\to \e_5 z$. The order of the full
symmetry group $\Z_5\times\Z_2$ equals 10.
This curve does not belong to any 
family of symmetric curves mentioned above.
\begin{figure}[htb]
    \centering 
     \includegraphics[width=10cm]{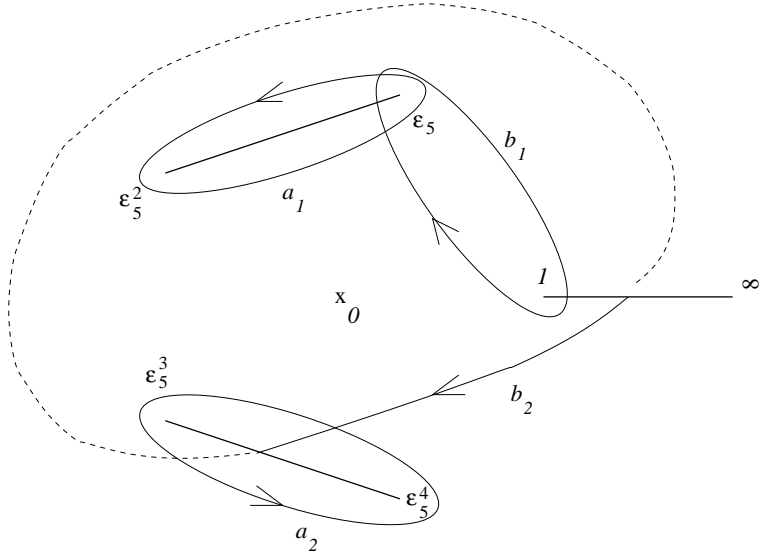}
    \caption{Canonical basis of cycles on the curve with $\Z_5\times\Z_2$ symmetry}
    \label{Z5fig}
\end{figure}
Choose the branch cuts and basic cycles $(a_i, b_i)$ on $\L$ as shown 
in Fig.~\ref{Z5fig}.  The generator $\mu$ acts on this basis as follows:
\be
b_1^{\mu}=-b_1+b_2\;, \hskip0.5cm b_2^{\mu}=-b_1+a_1+a_2\;, \hskip0.5cm
a_1^{\mu}=a_2\;, \hskip0.5cm a_2^{\mu}=-b_1\;,
\ee
therefore the action of the generator $\mu$ on
the vector $(b_1,\,b_2,\,a_1,\,a_2)^t$ is given by the following $\mod$ matrix:
\be
T^{\mu}=\left(\ba{cccc}   -1 & 1 & 0 & 0 \\
                          -1 & 0 & 1 & 1 \\
                           0 & 0 &  0 & 1 \\
                          -1 & 0 & 0 & 0 \ea \right)
\la{stabB3}
\ee
Invariance of the matrix of $b$-periods under the transformation 
defined by $T^{\mu}$  implies 
together with the positive definiteness of the matrix $\Im\B$ the following structure of 
this matrix in the chosen basis:
\be
\B=\left(\ba{cc} \e_5                & \e_5/(1+\e_5)  \\
                  \e_5/(1+\e_5)    & 1-\e_5^4 \ea\right)\;.
\la{Z5mat} 
\ee

The curves  (\ref{D6}), (\ref{Bolza}), (\ref{Z5})
are curves with large automorphism groups, i.e., they do not admit small deformations
preserving their groups of symmetries. They play an important role in the
subsequent analysis.

\section{$\F$ on the full moduli space of genus two curves}

\subsection{Boundary behavior and symmetry of $\F$}

Consider $\F$ (\ref{F}) as a function on Gottschling's fundamental domain 
$\G$. It is straightforward to prove 
the following lemma:
\begin{lemma}\la{lembou}
The function $\F$ (\ref{F})
vanishes if $\B\in\D$, as well as if $\B$ tends to the infinity of $\G$.
\end{lemma}
{\it Proof.} (see, for example \cite{Silhol}; we reproduce the proof 
here for completeness) The definition of the function of zero argument with vectors of characteristics
${\bf p,q}$, where ${\bf p}= (p_1,\,p_2)$ and   ${\bf q}= 
(q_1,\,q_2)$, reads:
\be
\Theta[^{\bf p}_{\bf q}](0|\B)=\sum_{{\bf m}\in \Z^2}
\exp\{\pi i \langle \B({\bf m}+{\bf p}),{\bf m}+{\bf p} \rangle + 
2\pi i\langle {\bf m}+{\bf p}, {\bf q}\rangle\}
\la{deftheta}
\ee
Suppose that $\B\in\D$; taking into account the  modular invariance of $\F$,
we can assume  $\B_{12}= 0$, while $\B_{11}$ and $\B_{22}$ remain finite.
Then the genus two theta function splits into the product of two 
genus one theta functions
with moduli $\B_{11}$ and $\B_{22}$: 
\be
\Theta[^{\bf p}_{\bf q}](0|\B)=\th[^{p_1}_{q_1}](0|\B_{11})\; \th[^{p_2}_{q_2}](0|\B_{22})
\ee
Therefore the even genus two theta constant corresponding to 
${\bf p}={\bf q}=(1/2,1/2)$ vanishes in this limit 
($\th\left[
\begin{smallmatrix}
    1/2  \\
    1/2
\end{smallmatrix}\right]
(0|\B_{11})=0$ since $\left[
\begin{smallmatrix}
    1/2  \\
    1/2
\end{smallmatrix}\right]$ is the odd genus one characteristic).

Now consider the part of the  boundary of $\M$, where the matrix of 
$b$-periods diverges, i.e., where some 
or all $y_i$ from (\ref{eq:num1}) tend to $+\infty$. Then ${\rm det}\Im\B$ diverges 
as a polynomial (of degree 2 with respect to $y_2$ and degree 1 with respect to $y_1$ and $y_3$),
while some theta constants vanish exponentially.

$\square$

To simplify the numerical analysis of the behavior of the function $\F$ on $\G$
we shall take into account the symmetry provided by the following lemma: 
\begin{lemma}
The function $\F$ has the following symmetry property:
\be
\F(-\overline{\B})=\F(\B)
\la{symB}\ee
\end{lemma}
{\it Proof.} This is a straightforward analog of the symmetry 
(\ref{mirgen1}) in the genus one case.
The proof is simple: each theta constant satisfies the relation $\Theta[\beta_s](-\overline{\B})=
\overline{\Theta[\beta_s](\B)}$ as a corollary of the definition (\ref{deftheta}); together with an
obvious symmetry of ${\rm det}\Im\B$ this implies (\ref{symB}).

$\square$

\subsection{Curves with large automorphism groups as critical points of $\F$}

Here we  show that all three Riemann surfaces ($D_6$, Burnside and $\Z_5$)  with large groups 
of automorphisms are critical points of any $\mod$ invariant function 
on $\S$ which is  real-analytic in open neighborhoods of the matrices of $b$-periods of these curves. 
In particular these Riemann surfaces are critical points 
of our function $\F$ (\ref{F}).

Consider an arbitrary (real-)analytic function $\F$ on $\S\setminus\D$. 
To analyze the power series of the function $\F$ in a neighborhood of 
some point of $\S$, 
it is convenient to map $\S$ to the generalized unit ball
$$\frak{U}=\{w\in M(2, {\mathbb C})\ :\ w=w^t; \ ww^*\leq I\}\;$$
in a way that the center of the power series is mapped to the origin.
Consider the Cayley transformation
$$z\mapsto w=(z-iI)(z+iI)^{-1}\;,$$
where $I$ is the $2\times 2$ unit matrix,
which gives a biholomorphic map $\S\rightarrow \frak{U}$.

To map a given point of $\S$ to the origin of  $\frak{U}$ we shall superpose the Cayley transformation 
with some automorphism of   $\frak{U}$.
Recall (see \cite{Siegel}) that all holomorphic automorphisms of the generalized 
unit ball $\frak{U}$ are given by the transformations
\be
w\mapsto (Aw+B)(Cw+D)^{-1},
\la{holaut}
\ee
where the $(2\times 2)$-matrices $A, B, C, D$ satisfy
the following constraints:
\begin{equation}\label{constr}
C=\overline{B},\ \ D=\overline{A}, \ \ \overline{A}A^t-\overline{B}B^t=I,\ \ AB^t=BA^t.
\end{equation}

A holomorphic  transformation of $\frak{U}$  to itself  such that a given point 
$S\in \frak{U}$ is mapped to the origin, looks as follows
 (see \cite{Siegel}, p. 177). Introduce an auxiliary matrix 
 \begin{equation}\label{koren0}
R=\sqrt{I_2-S\overline{S}},
\ee
and define matrices  $A, B, C, D$ by
\be\label{koren}
\ \ A=\overline{R^{-1}},\ \
B=-AS,\ \  C=\overline{B},\ \ D=\overline{A}.
\end{equation}
These  matrices  satisfy (\ref{constr}) and
$$(AS+B)(CS+D)^{-1}=0\;.$$

Consider a point $z_0$ of the Siegel half-space  $\S$ and map it to the point  $S=(z_0-iI_2)(z_0+iI_2)^{-1}\in\U$.
Define matrices $A_1, B_1, C_1, D_1$ 
by the equality
\begin{equation}
\begin{pmatrix} A_1 & B_1 \\ C_1 & D_1
\end{pmatrix}=
\begin{pmatrix} A & B \\ C  & D
\end{pmatrix}
\begin{pmatrix} I_2 & -iI_2 \\ I_2 & iI_2
\end{pmatrix}
\end{equation}
with $A, B, C, D$ from (\ref{koren}). 

Then 
the biholomorphic transformation ${\cal K}_{z_0}$ between the upper Siegel half-space and the generalized unit ball ${\cal K}_{z_0}:\S\rightarrow \frak{U}$
defined by
$${\cal K}_{z_0}(z)=(A_1z+B_1)(C_1z+D_1)^{-1},$$
maps the point $z_0\in \S$ to $0\in\U$.

\begin{theorem}\label{osnlem} Let $\F$ be any $\mod$-invariant function 
on $\S$, which is $C^1$ in open neighborhoods of the points 
$${\B}_1=\begin{pmatrix} \frac{i}{\sqrt{2}}-\frac{1}{2} & \frac{1}{2} \\ 
\frac{1}{2} & \frac{i}{\sqrt{2}}-\frac{1}{2}
\end{pmatrix}, \   \
{\B}_2=\frac{i}{\sqrt{3}}\begin{pmatrix} 2 & 1 \\ 1 & 2
\end{pmatrix} \  {\rm and}\
{\B}_3=\begin{pmatrix} \epsilon_5 & \frac{\epsilon_5}{1+\epsilon_5} \\ \frac{\epsilon_5}{1+\epsilon_5} & 1-\epsilon_5^4
\end{pmatrix},$$
which are the matrices of $b$-periods of the curves with large automorphism groups $y=x(x^4-1)$, $y=x^6-1$ and $y=z^5-1$ respectively.
Then ${\rm grad}\,\F$ vanishes at the points $\B_1$, $\B_2$ and $\B_3$.
\end{theorem}
{\bf Proof.}

\begin{enumerate}
\item
{\bf Point $\B_1$}. Let $T:=T^{\mu_1}$ where $ T^{\mu_1}$ is the 
first of the matrices  (\ref{stabB1})
which leave the point $\B_1$ invariant.

The function $g\ :\ \frak{U}\rightarrow {\mathbb C}$ defined by
\begin{equation}\label{fung} 
g(w):=\F({\cal K}_{\B_1}^{-1}w)
\end{equation}
satisfies the equation
\begin{equation}\label{inv}
g(w):=g({\cal K}_{\B_1}T {\cal K}_{\B_1}^{-1}w)\, .
\end{equation}
One can check that
\begin{equation}\label{unitar}
{\cal K}_{\B_1}T {\cal K}_{\B_1}^{-1}w=U_1 w U_1^{t}
\end{equation}
with a unitary matrix $U_1$ (recall that a $4\times 4$ matrix ${\cal K}_{\B_1}T {\cal K}_{\B_1}^{-1}$ having the block structure
$$\begin{pmatrix}A & B\\ C& D \end{pmatrix}
$$
acts on a $2\times 2$ matrix $w$ according to (\ref{holaut})).
Let
$$w=\begin{pmatrix} x & z\\ z & y \end{pmatrix},
\  \ U_1=\begin{pmatrix} u & v \\ w & t\end{pmatrix} \ {\rm and} \ \
U_1 w U_1^t=\begin{pmatrix} \tilde{x} & \tilde{z}\\ \tilde{z} & \tilde{y} \end{pmatrix}\;;
$$
 then
$$\begin{pmatrix}\tilde{x}\\ \tilde{y}\\ \tilde{z}\end{pmatrix}=
{\mathbb A}\begin{pmatrix}x\\ y \\z \end{pmatrix}$$
with
\begin{equation}\label{TRANSFORM}
{\mathbb A}=\begin{pmatrix}
u^2 & v^2 & 2uv\\
w^2 & t^2 & 2tw\\
uw & vt & vw+ut
\end{pmatrix}\, .
\end{equation}
The eigenvalues of the matrix ${\mathbb A}$ are given by $\{-1, i, -i\}$.

Due to (\ref{inv}) the partial derivatives  $g_x$, $g_y$ and $g_z$ of the function $g$ at
the point $w=0$ satisfy the equation
$$\begin{pmatrix}g_x\\ g_y \\g_z\end{pmatrix}={\mathbb A}^t\begin{pmatrix}g_x\\ g_y\\ g_z\end{pmatrix}.$$
Since $1$ is not an 
eigenvalue of the matrix ${\mathbb A}$, all these derivatives vanish 
for the function $\F$.

\item {\bf Point ${\B}_2$.} 
Here the same scheme applies to the matrix $T^{\mu_1}$ from 
(\ref{stabB2}), which belongs to the stabilizer of the point $\B_2$. 
In this case the spectrum  of the corresponding matrix ${\mathbb A}$ looks as follows:
$${\rm spectrum}\,({\mathbb A})=\left\{-1,\e_3,\e_3^2\right\}\;;$$
since again $1$ is not in the spectrum, ${\rm grad}\,g$ vanishes at the origin and ${\rm grad}\F$
vanishes at $\B_2$.

\item {\bf Point ${\B}_3$.} In this case we use the  modular 
transformation (\ref{stabB3}) which leaves the point $\B_3$ 
invariant. Then the spectrum of the corresponding matrix ${\mathbb 
A}$ reads:
$${\rm spectrum}\,({\mathbb A})=(\e_5,\e_5^2,\e_5^4)
\;;$$
since again $1$ is not among the eigenvalues, 
 ${\rm grad}\,g$ vanishes at the origin. Thus ${\rm grad}\F$ vanishes at $\B_3$.
\end{enumerate}

$\square$

\subsection{Numerical results}

The rigorous results established above can be summarized as follows: 
the function $\F$
 vanishes at the boundary of $\M$ and has at least
three critical points  corresponding to  Riemann surfaces with large automorphism groups. 
As we shall see below, knowing the value of the virtual Euler characteristic of $\M$, 
one can prove the existence of another critical point of $\F$.

However, since it was only possible to 
completely characterize the extremal points of even genus 1 modular functions numerically 
(see section 2), there is little hope for the time being to obtain 
complete analytical results for the decisively more complicated case of genus 2. Thus 
we turn to a numerical approach  to search for critical points of $\F$ in $\G$. 

The unboundedness of the domain $\G$ in the directions
of $y_1,\,y_2$ and $y_3$ (\ref{eq:num1}) is not a problem since 
the function $\F$ decreases exponentially for large $y_{i}$.  It 
turns out that a restriction to value of $y_{i}$ with
$y_i\leq2$ within the Gottschling 
domain is sufficient (notice that the algorithm explained below finds 
critical points with bigger values of the $y_{i}$, but these do not 
lie in the fundamental domain).
Furthermore we use  the symmetry (\ref{symB}) which in terms of 
$\{x_i,y_i\}$ looks as follows: $\F(\{-x_i,y_i\})=\F(\{x_i,y_i\})$ to 
decrease the amount of 
computation by a factor of $2$: in addition to Gottschling's conditions we assume $x_3>0$.

To locate the critical points of $\F$ inside of $\G$ we need to analyze the 
length of ${\rm grad}\F$ in $\G$. For the numerical evaluation of  ${\rm grad}\F$ at a given point of 
$\M$, we first differentiate $\F$ analytically with respect to the $\B_{ij}$. The differentiation of 
${\rm det}\Im \B$ is obvious; to differentiate theta constants  we  differentiate term by term
the series (\ref{deftheta}) to get:
$$
\f{\p}{\p\B_{jk}}\Theta[^{{\bf p}}_{{\bf q}}](0|\B)
$$
$$=\sum_{m_1,m_2\in\Z}
2\pi i (m_j+p_j)(m_k+p_k)/(1+\delta_{jk})
\exp\left\{\pi i\sum_{l,n=1,2}\B_{ln}(m_l+p_l)(m_n+p_n) +2\pi i 
\sum_{l=1,2}(m_l+p_l)q_l\right\}\;.
$$

The theta functions and its derivatives were approximated numerically 
via finite sums, $|m_{1,2}|\leq N$, in Matlab. For Riemann matrices in the 
fundamental domain, values of $N=3,4$ were sufficient to reach 
machine precision\footnote{Matlab works with a precision of 16 
digits; due to rounding errors machine precision is typically limited 
to 14 digits.}. We covered Gottschling's domain with a cartesian 
grid with 40 points in each direction (since we restricted the 
analysis to positive values of $x_{3}$ only 20 points for 
this direction were needed in the computation to obtain the same resolution as for 
the other $x$-directions). It turns out that roughly 40 \% of the 
points lie inside the fundamental domain.
The modular invariants were only calculated at these points.
In a first step we 
numerically identify values close to the minimum of the gradient 
for a given value of 
$y_{1}$ (within 0.01 of the respective minimum, in total more than 
8000 terms). The found values are used as an initial guess to 
search for a 
zero of the gradient. To identify the stationary points we use the 
algorithm of \cite{optim} which is implemented as the function 
\emph{fminsearch} in Matlab.

It turns out that four out of six critical points found numerically 
coincide to the order of machine precision with the points listed in Theorem \ref{osnlem} 
(which provides an additional test of the numerics).
Two other points are $\mod$ equivalent and  correspond to some 
Riemann surface from the $D_3$ family. All these points are located 
on the boundary of the fundamental domain.

The full list of the critical points found in one half of Gottschling's
domain, where $x_3\geq 0$, reads: 
\begin{enumerate}
\item
A point which can be identified with
\begin{equation}
    \begin{pmatrix}
        \eta & (\eta-1)/2  \\
        (\eta-1)/2 & \eta
    \end{pmatrix},
    \label{burnextr}
\end{equation}
where $\eta=(1+2\sqrt{2}i)/3$; this point is $\mod$ equivalent to
the point $\B_1$ of Theorem \ref{osnlem}, the Burnside curve.
The value of $\F$ for this point equals
$0.3106$ (we only give here 4 digits for the sake of presentation 
though at least 13 digits are known); this is the global maximum of 
the function $\F$. This result is confirmed by the
computation of  the signature of the Hessian which equals $(0,6)$.

\item
A point which can be identified with the point $\B_2$ of Theorem \ref{osnlem},
i.e., with the matrix of $b$-periods of the $D_6$ curve. The value of $\F$ at this point equals  
$0.2507$;  the signature of the Hessian is $(3,3)$.

\item
Two points which can be identified with the
matrices of $b$-periods of the $\Z_5$ curve: 
\begin{equation}
    \begin{split}
    \B=&
    \begin{pmatrix}
        \e_5 & \e_5+\e_5^{3}  \\
        \e_5+\e_5^{3} & -\e_5^4
    \end{pmatrix},\quad 
    \begin{pmatrix}
        -\e_5^4 & \e_5+\e_5^{3}  \\
        \e_5+\e_5^{3} & \e_5    
        \end{pmatrix}
\end{split}
    \label{Z5extr},
\end{equation}
where $\e_5=\exp(2\pi i/5)$. These two points are $\mod$ equivalent to the point $\B_3$
of Theorem \ref{osnlem}.
The value of the function $\F$ at these points equals $0.2912$.

The modular equivalent points (\ref{Z5extr}) belong to the boundary 
of the fundamental domain,
and coincide when the boundary points get appropriately identified.
The signature of the Hessian at these points equals $(2,4)$. 

\item Two $\mod$ equivalent  points 
\begin{equation}
    \mathbf{B} = 
    \begin{pmatrix}
        1.0517i  & \pm 0.5 + 0.5259i \\
        \pm 0.5 + 0.5259i & 1.0517i
    \end{pmatrix}
    \label{D3extr},
\end{equation}
(we give the components of this matrix with higher precision 
below) which belong to the $D_3$ family.
The value of $\F$ at these points equals $0.3011$.
These points also belong to the boundary of $\G$.
 The signature of the Hessian equals $(1,5)$.
\end{enumerate}

The symmetry of $\F$ implies the existence of further critical points 
for negative values of $x_3$ in the fundamental domain except for the 
points above with $x_3=0$. However these additional points are again
related to their originals by $\mod$ transformations, since all of these points 
belong to the boundary of $\G$. Therefore this does not contradict the fact that
$\G$ is a fundamental domain of $\mod$.

\begin{remark}\rm
The existence  of the  critical point (\ref{D3extr})  of $\F$ might 
appear surprising  since it
does not follow from Theorem \ref{osnlem}: this is not a Riemann 
surface with a large 
automorphism group. This phenomenon is also new in comparison with genus one, where
both critical points
of $\fo$ (\ref{Fgen1}) correspond to curves with large automorphism groups. However
the existence of such a point in genus two 
follows from the mass formula for the virtual (orbifold) Euler characteristic 
of $\M$ (see \cite{Schmutz0,Schmutz}  and section \ref{Euler} below). 
In \cite{Schmutz}, too, a point from the $D_3$ family was identified 
as a critical of the function $syst$. We do not know whether the critical $D_3$ curve
from \cite{Schmutz} coincides with (\ref{D3extr}). In our attempt to understand whether the point
(\ref{D3extr}) has a universal character we considered another $\mod$ invariant function on $\S$ - the 
absolute value of the first Igusa invariant \cite{Igusa}. Numerical analysis shows that 
the curve  (\ref{D3extr}) is no longer critical for this new function on $\M$. 
Therefore it remains an interesting problem to find some geometrical interpretation of 
the curve (\ref{D3extr}) (for instance the analogous curve of 
\cite{Schmutz}, which is critical for $syst$,
was claimed to be arithmetic).
The value of $r$ in equation (\ref{D3}) defining the extremal $D_3$ 
curve can be computed by
using Rosenhain's formulas (see e.g.~\cite{buser}) 
for branch points in terms of the matrix of $b$-periods; this value is given by
(all digits given below are reliable, except for the last two):
$$
r=0.22373907612077
$$
\end{remark}


\section{$\F$ on different strata of $\M$}

\subsection{Curves with  $\Z_2$ symmetry}

For the curves with $\Z_2$ reduced symmetry group (\ref{Z2}), when the matrix of $b$-periods 
for the canonical basis of cycles shown in Fig.  \ref{Z2fig},
has the form of (\ref{z2prym}), it is possible to express the function $\F$ 
(\ref{F}) in terms of elliptic theta functions of the moduli $\a$ and 
$\b$ by using the reduction formula
for genus two theta functions.
The result is given by the following lemma:

\begin{lemma}
Let the matrix of $b$-periods of a genus two curve have the form (\ref{z2prym}). Then
the function $\F$ (\ref{F}) can be represented in terms of genus one 
theta functions as
follows:
\begin{equation}    
\F =  \frac{1}{4} \fo^3(\a)\fo^3(\b)
(\Im \a)(\Im \b)\left|(\th_{3}^{4}(\a)  \th_{4}^{4}(\b)
-\th_{4}^{4}(\a)\th_{3}^{4}(\b)\right|\label{FZ2}.
\end{equation}
where $\fo(\a)$ is given by (\ref{Fgen1}); $\th_i(\a)$, $i=2,3,4$ are 
the genus one theta constants
of module $\a$.
\end{lemma}
{\it Proof.} When the matrix of $B$-periods has the form (\ref{z2prym}), the genus two
theta function decomposes into a combination of elliptic theta functions with
moduli $2\a$ and $2\b$: (the first argument of
the theta functions below is 0)
\begin{equation}    
\Theta  \begin{bmatrix}     a & b  \\       c
& d     \end{bmatrix}(\B)=\sum_{e \in\{0,1/2\}}^{}\vartheta
\begin{bmatrix}             (a+b)/2+e   \\                  c +d
\end{bmatrix}\vartheta(2\b)     
\begin{bmatrix}             (a-b)/2+e
\\                  c -d                
\end{bmatrix}(2\a)
\label{eq:thetaprym} 
\end{equation}
By regrouping the ten theta constants from (\ref{F}) into five pairs, 
and by  using
the inverse binary addition formula for elliptic theta functions:
\begin{equation}   
 \vartheta      \begin{bmatrix}     a
   \\       b   \end{bmatrix}(2\a)\vartheta     \begin{bmatrix}     c
   \\       b   \end{bmatrix}(2\a)=\frac{1}{2}\sum_{d\in\{0,1/2\}}^{}
   \exp(-i\pi ad)\vartheta      \begin{bmatrix}     a+c  \\
   b+d  \end{bmatrix}(\a)\vartheta      \begin{bmatrix}     a-c  \\
   d    \end{bmatrix}(\a)   
\label{inverse},
\end{equation}
we get (\ref{FZ2}) after elementary manipulations.

$\square$


It is clear that the Burnside curve is also the absolute maximum of 
$\F$ on $\M(\Z_2)$; the 
$D_6$ curve and the extremal $D_3$ curve must also be  critical points of $\F$ on $\M(\Z_2)$
(the $\Z_5$ curve does not belong to $\M(\Z_2)$).
The question whether $\F$ has critical points on $\M(\Z_2)$ in 
addition to these three, 
can also be given only a numerical (but rather conclusive) answer.

Taking into account  Theorem \ref{funz2} and Lemma \ref{simfunz2}, we study
numerically the behavior of $\F$ and ${\rm grad}\F$ when  $(\a,\b)\in\Omega\times\Omega(2)$.
The numerical analysis with 100  points in each direction and a 
subsequent refined analysis as described above identifies the following 
critical points on  $\Omega\times\Omega(2)$, all of which coincide 
up to $\mod$ transformations with one of the four 
critical points of $\F$ on the whole $\M$:
\begin{enumerate}
\item
The point which coincides with the matrix of  $b$-periods of 
the Burnside curve:
\begin{equation}
 \a=\sqrt{2}i \;\hskip1.0cm
\b=\f{2}{3}+\sqrt{2}i  
    \label{burnz2}.
\end{equation}
This point, of course, gives the global maximum of $\F$ on $\M(\Z_2)$.
The signature of the Hessian at this point is $(0,4)$.

\item
Five points corresponding to the $D_6$ curve. The first 
point  coincides with the point $\B_2$ of Theorem \ref{osnlem}.
The following four  points which also belong to $\Omega\times\Omega(2)$ are 
$\mod$ equivalent to  $\B_2$
(i.e., they also represent the $D_6$ curve):
$$
\a=\pm \f{1}{2}+\f{i\sqrt{3}}{2} \hskip0.8cm   
\b=\mp \f{1}{2}+\f{i\sqrt{3}}{2}
$$
\be
\a=\pm \f{3}{2}+\f{i\sqrt{3}}{2} \hskip0.8cm   
\b=\mp \f{1}{2}+\f{i\sqrt{3}}{2}
\la{D6z2}
\ee
$$
\a=\pm \f{3}{2}+\f{i\sqrt{3}}{2} \hskip0.8cm   
\b=\mp \f{1}{2}+\f{i\sqrt{3}}{2}
$$
The signature of the Hessian at these points is $(2,2)$.


\item
Three points which
are $\mod$ equivalent to the critical $D_3$ curve (\ref{D3extr})
found on the full $\M$:
\begin{equation}
    \begin{split}
    &\begin{pmatrix}
        0.3835 + 0.7874i &-0.4339 + 
        0.2114i\\
         -0.4339 + 0.2114i&  0.3835 + 0.7874i
      \end{pmatrix},\\
      &\begin{pmatrix}
       1.0517i& -0.5 + 0.5259i\\
      -0.5 + 0.5259i&  1.0517i
     \end{pmatrix},\\
     &\begin{pmatrix}
         0.5+ 1.0517i&  0.5259i\\
         0.5259i & 0.5+1.0517i
        \end{pmatrix}.
      \end{split}
    \label{D3z2}
\end{equation}
The signature of the Hessian at these  points  is $(1,3)$.
\end{enumerate}

\begin{remark}\rm
We see that the fundamental domain $\Omega\times\Omega(2)$ is indeed bigger than the moduli space
$\M(\Z_2)$: there are  equivalences between different points of 
$\Omega\times\Omega(2)$ given by $\mod$ transformations which do not correspond to any 
$\Gamma\times\Gamma(2)$ transformation on $\Omega\times\Omega(2)$. 
\end{remark}

\subsection{Curves with $D_2$ symmetry}

The matrices of $b$-periods  (\ref{BD22}) of the $D_2$ curves form a subfamily 
of the two-parametric family (\ref{z2prym})
with $\b=\a-1$; since
$\vartheta_{3,4}(\a-1)=\vartheta_{4,3}(\a)$, we have in this case
\be
\F =  \frac{1}{4} \fo^6(\a)
(\Im \a)^2\left|(\th_{4}^{8}(\a)
-\th_{3}^{8}(\a)\right|
\label{FD2};
\ee
due to Lemma \ref{lemg2}, the moduli space $\M(D_2)$ in the variable $\a$
coincides with the
fundamental domain $\Omega_0(2)+$ of the group $\Gamma_0(2)+$.
  
\begin{figure}[htb]
    \centering 
    \includegraphics[width=10cm]{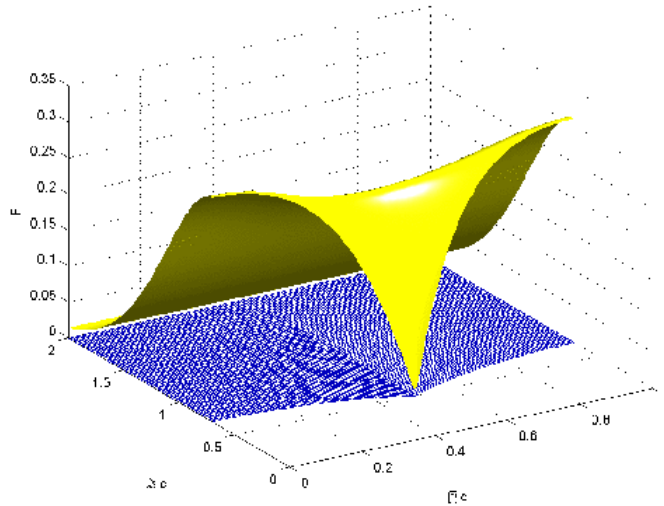}
    \caption{$\F(\s)$ in the fundamental domain $\Omega_0(2)+$}
   \label{plotD2}
\end{figure}

The plot of $\F$ in $\Omega_0(2)+$ is shown in Fig.~\ref{plotD2}.
It is clear that the points corresponding to the Burnside and to the $D_6$ curves, which are critical 
points of $\F$ from a six-dimensional point of view, are also critical points of $\F$
on $\M(D_2)$ from a two-dimensional point of view. 
No other critical points of $\F$ in $\M(D_2)$  were found.

Namely for $\a\in \Omega_0(2)+$ we find the following critical points 
(see Fig.~\ref{Omega200fig}):
\begin{enumerate}

\item Two points corresponding to the Burnside curve
are found at $\a=i/\sqrt{2}$ and $\a=i/\sqrt{2}+1$; the first of these
values gives the 
point $\B_1$ of Theorem \ref{osnlem}. 
This is the absolute maximum of $\F$, and the signature of the Hessian is $(0,2)$.

\item The point corresponding to the $D_6$ curve is given by 
the  value $\a= \frac{1}{2}+i\frac{\sqrt{3}}{2}$. The signature of the Hessian  at this point is $(1,1)$.
\end{enumerate}

The point of intersection of the two circles limiting the fundamental domain is
$\a=1/2+i/2$; at this point the genus two Riemann surface splits into
two tori; the function $\F$ vanishes at this point, as well as at $\a\to i\infty$.

\subsection{Curves with $D_3$ symmetry}

The matrices of $b$-periods  (\ref{BD3}) of this family form a subfamily 
of the two-parametric family (\ref{z2prym})
with $\a=3\s$ and $\b=\s$; then
\be
\F(\s) =  \frac{3}{4} \fo^3(\s)\fo^3(3\s)
(\Im \s)^2\left|(\th_{4}^{4}(\s)\th_{3}^{4}(3\s)
-\th_{3}^{4}(\s)\th_{4}^{4}(3\s)\right|
\label{FD3}\;;
\ee
due to Lemma \ref{lemg3}  the moduli space $\M(D_3)$ in the variable $\s$
coincides with the
fundamental domain $\Omega_0(3)+$ of the group $\Gamma_0(3)+$.

\begin{figure}[htb]
   \centering 
   \includegraphics[width=10cm]{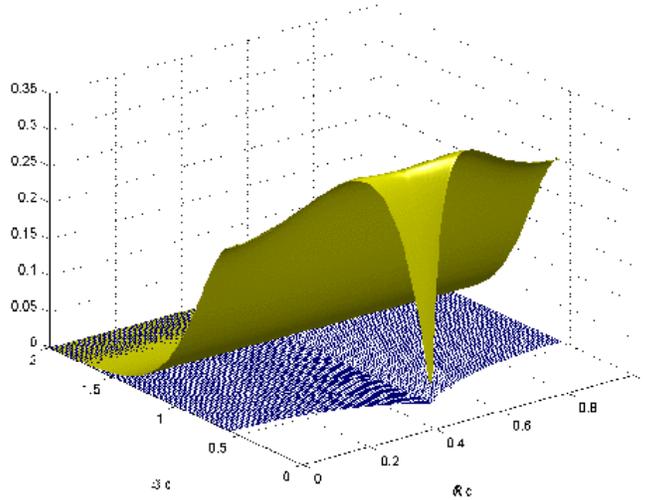}
   \caption{Plot of $\F$ in the fundamental domain $\Omega_0(3)+$}
   \label{plotD3}
\end{figure}

The plot of $\F$ in $\Omega_0(3)+$ is shown in Fig.~\ref{plotD3}.
It is clear that the points corresponding to the Burnside, the $D_6$ and the critical $D_3$ 
curves, which are critical 
points of $\F$ from a six-dimensional point of view, are also critical functions of $\F$
on $\M(D_3)$ from a two-dimensional point of view. 
No other critical points (i.e., such that the two-dimensional gradient vanishes) 
of $\F$ in   $\M(D_3)$  were found numerically.

Namely we find for $\s\in \Omega_0(3)+$ the following critical points 
(see Fig.~\ref{Omega300fig}):
\begin{enumerate}
\item
Two points corresponding to the Burnside curve are given by
\begin{equation}
    \s = \f{1}{3}+\f{i\sqrt{2}}{3}\;, \hskip0.8cm \s = \f{2}{3}+\f{i\sqrt{2}}{3}
    \label{D3burn}.
\end{equation}
the corresponding matrices of $b$-periods are related by $\mod$ transformations 
to the point $\B_1$ of Theorem \ref{osnlem}.
The signature of the Hessian is $(0,2)$.

\item
Two  points corresponding to the $D_6$ curve:
 $\s = i/\sqrt{3}$ and  $\s = 1+i/\sqrt{3}$; the first of these values
 gives  the point $\B_2$ of Theorem \ref{osnlem}.
The signature of the Hessian is $(1,1)$.

\item
The point corresponding to the extremal $D_3$ curve is given by:
\begin{equation}
    \s = 0.5+0.5259i.
    \label{D3D3}
\end{equation}
The corresponding matrices of $b$-periods are $\mod$ equivalent to the matrices of $b$-periods 
(\ref{D3extr}). 
The signature  of the Hessian is $(1,1)$.
\end{enumerate}

The point of intersection of the two circles limiting the fundamental domain
 $\s=\f{1}{2}+\f{i}{2\sqrt{3}}$ is the boundary point of the moduli
 space, where the genus two surface splits into two tori; according to
 Lemma \ref{lembou}, $\F$ vanishes there.

\section{$\F$ and  Euler characteristics of moduli spaces}\label{Euler}

Given a function on the moduli space with non-degenerate critical 
points, it is natural to ask whether 
one can use this function as a  Morse function, i.e., whether one can extract topological information
(say, the Euler characteristic, which we shall discuss here) about the
space. Since
$\M$  (as well as its symmetric strata) is an orbifold, 
one can define many different Euler characteristics of $\M$ (see, for example, \cite{Hirz,Bryan}).
We shall speak here about the ordinary Euler characteristic and the
orbifold Euler characteristic in the sense of \cite{Thurston} defined as follows.

Let $X$ be a manifold (perhaps, with non-trivial boundary  $\partial X$); let ${\cal G}$ 
be a group acting on $X$ effectively and properly discontinuously. Consider an arbitrary proper
 cell division of the orbifold $\M=X/G$: $X=\cup_{j=1}^n X_j$ 
(the cell division is called proper if all points of a given cell $X_j$ have the same 
stabilizer $G_j$). Then the orbifold Euler characteristic of $\M$ is defined as follows:
\be
\chi_o(\M)=\sum_{j=1}^n \f{(-1)^{{\rm dim}\; X_j}}{\#(G_j)}
\la{defchio}
\ee

Suppose now ${\cal F}:\;X\to {\bf R}$ to be a function on $X$ vanishing on the boundary $\partial X$,
invariant under the action of the group ${\cal G}$ such that all critical points of $f$ are isolated and
non-degenerate.
The set of critical points of the function ${\cal F}$ must be invariant under the action of the group 
${\cal G}$;
let us assume that this set consists of a finite number $N$ of orbits; denote by 
$P_1,\dots, P_N$ the critical points  representing each orbit. 
Denote by $i_1,\dots,i_N$ the indices of the critical points  
$P_1,\dots, P_N$ (i.e., the
numbers of negative eigenvalues of the corresponding Hessians).  Denote by ${\rm Stab}(P_j)$ 
the stabilizer of the critical point $P_j$, i.e., the subgroup of ${\cal G}$ leaving 
the point $P_j$ invariant; 
the order of the subgroup is denoted by $\#({\rm Stab}(P_j))$.

Then the orbifold Euler characteristic can be expressed via the following ``mass formula'' 
(see \cite{Ash,Schmutz} and references therein):
\be
\chi_o(X/{\cal G})=\sum_{j=1}^N \f{(-1)^{i_j}}{\#({\rm Stab}(P_j))}
\la{stabil}
\ee
(We recall that if $f$ is a Morse function with a finite number of critical points on a {\it manifold},
then the ordinary Euler characteristic is equal to $\sum (-1)^{i_j}$, where $i_j$ are indices of
the critical points.

A moduli space $\M_g$ can be represented as the factor  ${\cal 
T}_g/{\cal G}_g$ of the Teichm\"uller 
space ${\cal T}_g$ by the 
mapping class group ${\cal G}_g$.   Let ${\cal F}$ be a $C^1$ function on  ${\cal T}_g$
invariant with respect to the mapping class group ${\cal G}_g$; let $\L_1,\dots,\L_N$ be 
Riemann surfaces representing the orbits of the critical points of ${\cal F}$ on ${\cal T}_g$.
Then for $g > 2$ ${\rm Stab}(\L_j)$ is simply the group ${\rm Aut}(\L_j)$ of (holomorphic) automorphisms of 
$\L_j$ and the mass formula (\ref{stabil}) looks as follows:
\be
\chi_o(\M_g)= \sum_{j=1}^N \f{(-1)^{i_j}}{\#({\rm Aut}(\L_j))}\;\hskip1.0cm g>2\;.
\la{stabil1}
\ee

For $g=1,2$ any curve from $\M_g$ possesses the hyperelliptic 
involution, and the right-hand side 
of (\ref{stabil1}) should be multiplied by $2$, i.e.,
\be
\chi_o(\M_{g})= 2\sum_{j=1}^N \f{(-1)^{i_j}}{\#({\rm Aut}(\L_j))}\;\hskip1.0cm g=1,2\;.
\la{stabil2}
\ee

\subsection{Genus one}

The moduli space $\M_1$ can be obtained by 
standard gluing of the boundary of the fundamental domain $\Omega$ 
(Fig.~\ref{Omegafig}); as a result we get
the sphere with one puncture (the asymptotic cylinder corresponding to $\s\to\infty$) and two orbifold
points corresponding to $\s=i$ (of index 2) and $\s=e^{\pi i/3}$ 
(of index 3). The ordinary Euler characteristic of this space is  $1$.

The orbifold Euler characteristic of $\M_1$ equals $-1/6$, as follows from considering an obvious cell division; this result can be 
easily reproduced using the ``mass formula'' (\ref{stabil2}). 
Namely, the function $f(\s)$ (\ref{Fgen1}) has two critical points: $\s=i$ and $\s=e^{2\pi i/3}$;
the order of the automorphism group of the torus with period $i$ is
$4$, and the order of the automorphism group of the torus with period $e^{2\pi i/3}$ is $6$;
thus the ``mass formula'' (\ref{stabil2}) gives $2(1/6-1/4)=-1/6$.

\subsection{Genus two }

\subsubsection{Full moduli space}

The ordinary Euler characteristic of the moduli space $\M$ in genus two as well as of 
its symmetric strata (Bolza subspaces) were discussed in a recent paper \cite{Getzler}.
In particular $\chi(\M)=1$. 


The orbifold Euler characteristic of $\M$ equals $-1/120$  \cite{HZ}. 
On the other hand, computing the right-hand side of (\ref{stabil2}) for our function $\F$, 
we get (taking into account the
indices  of ${\cal F}$ on the Burnside, the $D_6$, the $\Z_5$ and the critical $D_3$ curves are equal to
6, 3, 4 and 5 respectively, and the corresponding number of automorphisms which are
48, 24, 10 and 12 respectively):
\be
2\left(\f{1}{48}-\f{1}{24}+\f{1}{10}-\f{1}{12}\right)=-\f{1}{120}\;,
\ee
in agreement with \cite{HZ} (this computation coincides with the one done in \cite{Schmutz}
where a similar analysis was performed for the function $syst$).

\begin{theorem}\la{fourthcurve}
Let $\F$ be any  modular invariant smooth function on $\S$ vanishing on the 
boundary of the moduli space  $\M=(\S\setminus\D)/\mod$ and such that all of its 
critical points are non-degenerate. Then in addition to
critical points at the three Riemann surfaces with 
large automorphism groups, $\F$ must have at least one more critical 
point within the $D_3$ family
(\ref{D3}) or the $D_2$ family (\ref{D2}).
\end{theorem}
{\it Proof} is elementary: subtracting any  combination  of $1/48$, $1/24$ and $1/10$
from $-1/240$ with different signs, we get 8 different numbers. 
All of them have too large denominators to arise only from curves from outside of $D_2$ and $D_3$ families.

$\square$

 Theorem \ref{fourthcurve} does not exclude the existence of other critical points 
of the function $\F$ such that their total contribution to the mass formula (\ref{stabil2}) 
vanishes.
However our numerical analysis shows that $\F$ has indeed only four 
critical points on $\M$: three standard ones, and one more from the $D_3$ family, 
and all of them are non-degenerate.

\subsubsection{Symmetric strata}

Although we shall  speak here only about the genus two case, most of the
construction makes sense for an arbitrary genus.
Denote by $\Sigma^G$ the set  of Riemann surfaces
having the group $G$ as their group of automorphisms. The space $\Sigma^G$ contains also 
Riemann surfaces whose groups of automorphisms has $G$ as a subgroup. 
We shall denote by the same letter the element of the group and the biholomorphic automorphism of the Riemann surface corresponding to this element. 
All the homeomorphisms below are assumed to be orientation preserving.

Consider a basepoint $\L_0\in \Sigma^G$ and define a {\it $G$-symmetric marking} of a 
surface $\L\in\Sigma^G$ as a homeomorphism
$$\phi\;:\;\L_0\rightarrow \L$$
such that for any element $g_1$ of $G$ 
the map 
 $\phi g_1\phi^{-1}$ is homotopic  to some  $g_2\in G$.
Two $G$-symmetrically marked surfaces $\L_1$ and $\L_2$ from $\Sigma^G$ with $G$-symmetric markings 
$\phi_1: \L_0\rightarrow \L_1$ and $\phi_2 : \L_0\rightarrow \L_2$ are said to be equivalent if
there exists a biholomorphic map $h : \L_1\rightarrow \L_2$ which is homotopic to the 
map $\phi_2\phi_1^{-1}$.

The set of equivalence classes of $G$-symmetrically marked surfaces from $\Sigma^G$ is called
{\it $G$-special Teichm\"uller space} and denoted by ${\cal T}^G$. We are not aware of any 
general results about the structure of these spaces; for genus two their description is given in
\cite{Schiller}.

We call a homeomorphism $\Phi:\;\L_0\rightarrow \L_0$  G-symmetric if for any $g_1\in G$ there 
exists $g_2\in G$ such that the homeomorphisms $\Phi g_1$ and
$g_2\Phi$ are homotopic.
Denote the group ${\cal G}_G'$ to be  the group of homotopy classes of all $G$-symmetric homeomorphisms
$\Phi:\L_0\rightarrow \L_0$. 

The group ${\cal G}_G'$ acts on ${\cal T}^G$ as follows.
Let $\L$ be a $G$-symmetrically marked surface with $G$-symmetric marking $\phi:\;\L_0\rightarrow \L$
and $[\Phi]$ be an element of ${\cal G}_G'$ defined by a $G$-symmetric homeomorphism 
$\Phi:\;\L\rightarrow \L$. Then $[\Phi](\L)$ is the same surface $\L$ marked as 
\begin{equation}\label{sm}
\phi \circ \Phi^{-1} :\;\L_0\rightarrow \L\;;
\end{equation}
obviously, this marking is also $G$-symmetric.
Let $I$ be a subgroup of ${\cal G}_G'$ whose elements preserve  any $\L\in {\cal T}^G$.

The G-{\it special mapping class group} ${\cal G}_G$ is defined as the factor-group 
${\cal G}_G={\cal G}_G'/I$ (this definition excludes the hyperelliptic involution in genus two).
The action of ${\cal G}_G'$ on ${\cal T}^G$ gives rise to 
the properly discontinuous effective action of 
${\cal G}_G$ on ${\cal T}^G$;  the orbifold ${\cal M}^G={\cal T}^G/{\cal G}_G$ is called 
G-{\it  special moduli space}.

Consider now the moduli spaces $\M(\Z_2)$ ($=\M^{\Z_2\times\Z_2}$), 
$\M(D_2)$ ($=\M^{D_2\times\Z_2}$) and $\M(D_3)$  ($=\M^{D_3\times\Z_2}$).
These moduli spaces are factors of the corresponding special Torelli spaces 
${\cal S}(\Z_2)$ (\ref{fd22}), ${\cal S}(D_2)$ (\ref{SD2}) and ${\cal S}(D_3)$ (\ref{SD3}) by the 
special Torelli groups (\ref{grG}), $\Gamma_0(2)+$ and $\Gamma_0(3)+$ respectively.

The universal covering of ${\cal S}(\Z_2)$ is the special Teichm\"uller space 
${\cal T}(\Z_2)=H\times H$ \cite{Schiller}. The special Teichm\"uller spaces 
${\cal T}(D_2)$ and ${\cal T}(D_3)$ are the
universal coverings of ${\cal S}(D_2)$ and ${\cal S}(D_3)$, respectively; they both  
coincide with the upper half-plane $H$. 

We denote corresponding special
mapping class groups by ${\cal G}(\Z_2)$, ${\cal G}(D_2)$ and ${\cal G}(D_3)$, respectively.
Then the moduli spaces $\M(\Z_2)$, $\M(D_2)$ and $\M(D_3)$ are represented as 
factors ${\cal T}(\Z_2)/{\cal G}(\Z_2)$,  ${\cal T}(D_2)/{\cal G}(D_2)$ and 
${\cal T}(D_3)/{\cal G}(D_3)$ respectively.


Let $\L$ be a surface from $\M^G$ and let $F$ be its {\it full} group of holomorphic automorphisms.
Then $G$ is a subgroup of $F$ and we denote by $N(G; F)$ the normalizer of $G$ in $F$. 
Choose some $G$-symmetric marking $\phi:\;\L_0\rightarrow \L$ of $\L$ and denote by 
${\rm Stab}(\L)$ the stabilizer of
the point $\L\in\M^G$ in the special mapping class group ${\cal G}_G$.

\begin{lemma}\la{STNORM}
The map $\kappa:\; N(G; F)\rightarrow {\rm Stab}(\L)$ defined by 
$$\kappa(f)=[\phi^{-1}f\phi]\;,$$
where $f\in N(G; F)$ and $[\dots ]$ denotes the homotopy class, 
is the homomorphism of two  groups.
\end{lemma}

{\it Proof.} Let $\Phi:=\phi^{-1} f\phi$.
To verify that $\kappa$ is the group homomorphism 
one has to show that 
\begin{enumerate}
\item  homeomorphism $\Phi$ is $G$-symmetric
\item 
 $[\Phi](\L)=\L$.
\end{enumerate}
Choose $g_1\in G$; then $\phi g_1\phi^{-1}=g_0$ for some $g_0\in G$. 
Since $f\in N(G; F)$, there exists some $g_2\in G$ such that $fg_0=g_2f$.
Therefore, $g_2f=fg_0\sim f \phi g_1 \phi^{-1}$ (where $\sim$ denotes 
the homotopy equivalence) and there exists a $g_3$ such that
$$\Phi g_1\sim \phi^{-1}g_2f\phi=\phi^{-1}g_2\phi\phi^{-1}f\phi\sim g_3\Phi\,,$$
which proves the $G$-symmetry of $\Phi$.

To prove the second statement we write
$(\phi \Phi^{-1})\phi^{-1}=\phi\phi^{-1}f^{-1}\phi\phi^{-1}=f^{-1}\;;$
thus the $G$-symmetric markings
$\phi:\;\L_0\rightarrow \L$ and $\phi\circ\Phi^{-1}:\;\L_0\rightarrow \L$ are equivalent. 
 
 $\Box$

\begin{lemma}
The groups $N(G;F)/G$ and ${\rm Stab}(\L)$ are isomorphic.
\end{lemma}
{\it Proof.} In view of the previous lemma it is sufficient to show that
\begin{enumerate}
\item $\kappa$ is a surjection,
\item ${\rm ker} \,\kappa=G$.
\end{enumerate}

1) Let $r\in {\rm Stab}(\L)$, then $r=[\Phi]$, where $[\Phi]$  is  a $G$-symmetric homeomorphism 
$\Phi: \L_0\rightarrow \L_0$ such that
\begin{equation}\label{e2}
(\phi\Phi^{-1})\phi^{-1}\sim h
\end{equation}
with some $h\in F$.
To prove the surjectivity of $\kappa$ we have to show that $h\in N(G; F)$ and $\kappa(h)=r$.
Let $g_1\in G$, then
$$g_1h\sim g_1 \phi \Phi^{-1}\phi^{-1}\sim\phi g_2\Phi^{-1}\phi^{-1}\sim\phi\Phi^{-1}g_3\phi^{-1}\sim
\phi\Phi^{-1}\phi^{-1}g_4\sim hg_4,$$ 
with some $g_2, g_3, g_4\in G$.
 Since two homotopic holomorphic automorphisms must coincide
(see e.g.  Lemma 6.5.5 from \cite{buser}), we have
$g_1h=hg_4$, and, therefore, $h\in N(G; F)$.
From (\ref{e2}) we get
$$r=[\Phi]=[\phi^{-1}h\phi]=\kappa (h)$$
as was stated.

2a)  $G\subset {\rm ker}\,\kappa$.

We shall prove that for any $g\in G$ the element $\kappa(g)$ of the special mapping class group
acts trivially on ${\cal T}^G$ and, therefore, is the unity (since the  special mapping class group action is effective).
Let $\L$ be a $G$-symmetrically marked surface with symmetric marking $\psi:\;\L_0\rightarrow \L$.
Then 
$$(\psi\phi^{-1}g^{-1}\phi)\psi^{-1}\sim g^{-1}$$
and, therefore, $\kappa(g)(\L)=\L$.

2b)  ${\rm ker}\,\kappa\subset G$.

Let $f\in {\rm ker}\,\kappa$ and $\L$ be a {\it generic} surface 
 admitting a holomorphic $G$-action (i.e. ${\rm Aut}(\L)=G$) and
let $\phi:\;\L_0\rightarrow \L$ be a $G$-symmetric marking.
Since $\kappa(f)={\bf 1}$, one has
$$(\psi\phi^{-1}f^{-1}\phi)\psi^{-1}\sim g$$
with some $g\in G$.
Therefore $f^{-1}\sim \phi \psi^{-1}g\psi\phi^{-1}\sim g_1$
with some $g_1\in G$
and Lemma 6.5.5 from \cite{buser} implies that $f^{-1}=g_1$ and, therefore, $f\in G$.

$\square$

\begin{theorem}
Let $\F$ be any smooth function on the special Teichm\"uller space ${\cal T}^G$ ($G$ is 
any of the groups $D_2\times\Z_2$, $D_3\times\Z_2$ or $\Z_2\times\Z_2$)
invariant with respect to the special mapping class group ${\cal G}_G$ and vanishing on the
boundary of ${\cal T}^G$. Suppose that the function $\F$ has critical points at
Riemann surfaces $\L_k$ which possess the automorphism groups $G_k$
respectively ($G$ is a subgroup of $G_k$). Denote by $H_k$ the factor-group of the normalizer $N(G,G_k)$ by $G$. Then  the orbifold Euler
characteristic of $\M^G$ is given by:
\be
\chi_o(\M^G)= \sum_k \f{(-1)^{i_k}}{\#(H_k)}
\la{defMG}
\ee
where $i_k$ is the index of the critical point $\L_k$ and $\#(H_k)$ is
the order of the subgroup $H_k$.
\end{theorem}

{\it Proof.} The mass formula (\ref{defMG}) immediately follows from  the mass formula (\ref{stabil})
and lemma \ref{STNORM}.

$\Box$
  
Let us now consider the three symmetric strata of $\M_2$ separately.


{\bf Space $\M(D_2)$.} 

Let us compute $\chi_o(\M(D_2))$ using the formula (\ref{defMG}).
Namely, let us enumerate the branch points $(0,\infty,1, i,-1,-i)$ of the Burnside curve $z(z^4-1)$ 
by the numbers $(1,2,3,4,5,6)$ respectively. Then the automorphism $\mu_1\;:\;z\to iz$ acts on the set
of branch points as the
permutation $(3,4,5,6)$. The involution $\mu_2\;:\; z\,\to\, (z+1)/(z-1)$ acts as $(15)(23)(46)$;
finally, the involution $z\to -1/z$ acts as $(12)(35)$.

The $D_2$ subgroup of $S_4$ is generated by the symmetries $z\to -z$ and $z\to r/z$ on the 
Riemann surface  $z(z^2-1)(z^2-r^2)$. On the  
 set of branch points $(0,\infty,1,r,-1,-r)$ these symmetries act by 
 the permutations $(35)(46)$
and $(12)(34)(56)$, respectively. The normalizer of the $D_2$ subgroup in $S_4$
contains one more 
generator: $(12)(46)$. This element generates the factor-group $H$ of the normalizer by $D_2$
which is isomorphic to  $\Z_2$.
Thus $\#(H)$ for the Burnside curve equals $2$.

Consider now the $D_6$ curve.
Let us  enumerate the branch points of the $D_6$ curve by $1,\dots, 6$
starting from $z=1$ counterclockwise. Then the involution $\mu_1\,:\;z\to\epsilon_6 z$
acts on this set as the permutation $(1,2,3,4,5,6)$; the involution $\mu_2\,:\;z\to -1/z$ acts as
the permutation $(14)(23)(56)$.  The $D_2$ subgroup of $D_6$ is generated by the
involutions $z\to -z$ and $z\to -1/z$ which correspond to the permutations 
$(14)(25)(36)$ and $(14)(23)(56)$. The normalizer of the $D_2$ subgroup in $D_6$ coincides 
with the subgroup itself, thus $\#(H)$ at $D_6$ curve equals $1$.

Thus, the formula (\ref{defMG}) gives 
$$\chi_o(\M(D_2))=-1+\f{1}{2}=-\f{1}{2}\;\;.$$

{\bf Space $\M(D_3)$.}

The $D_3$ symmetry on the curve $(z^3-1)(z^3-r^3)$ is generated by 
the transformations
$z\to r/z$ and $z\to \epsilon_3 z$. On the branch points of the Burnside curve these transformations act
as the permutations $(2,5)(4,6)(1,3)$ and $(3,2,6)(1,5,4)$, respectively. The normalizer of the
$D_3$ subgroup of $S_4$ coincides with the $D_3$ subgroup itself; 
thus at the $S_4$ curve $\#(H)=1$.
On the other hand, $D_3$ is the subgroup of index $2$ in $D_6$; thus it is a normal subgroup
and at the $D_6$ curve $\#(H)=2$.

The formula (\ref{defMG}) gives in this case
$$\chi_o(\M(D_3))=-1+1-\f{1}{2}=-\f{1}{2}\;.$$

{\bf Space $\M(\Z_2)$}. 

To compute $\chi_o(\M(\Z_2))$ via  (\ref{defMG})
we need to find the subgroups $H$ for the stationary points of $\F$, 
i.e., for 
the Burnside curve, the $D_6$ curve and the extremal $D_3$ curve.

The $\Z_2$ symmetry $z\to -z$ on the Burnside curve acts as the permutation $(35)(46)$ of the branch points;
using Maple we verify that the factor-group of the normalizer of this $\Z_2$ subgroup in $S_4$
by the $\Z_2$ subgroup is generated by the permutations  $(12)(34)(56)$ and $(12)(46)$   and has order 4.

The $\Z_2$  symmetry on the $D_6$ curve is $z\to \epsilon_6/z$; on the branch points it acts as
the permutation $(12)(36)(45)$. The factor-group of the normalizer of the $\Z_2$ subgroup by
the subgroup itself is another $\Z_2$ subgroup generated by the permutation $(14)(25)(36)$;
thus here $\#(H)=2$.

Consider a $D_3$ curve $(z^3-1)(z^3-r^3)$ and enumerate the branch points 
$(1,r,\epsilon_3,\epsilon_3 r, \epsilon_3^2,\epsilon_3^2 r)$ by the numbers $(1,\dots,6)$.
Then the $D_3$ group is generated by the symmetries $z\to \epsilon_3 z$ and $z\to r/z$;
the corresponding permutations of the branch points are $(135)(246)$ and $(12)(36)(45)$; the
permutation $(12)(36)(45)$ defines the $\Z_2$ subgroup. The normalizer of this subgroup coincides
with the subgroup itself; thus $\#(H)=1$.

Since the index of the Burnside curve in the $\Z_2$ family is $4$, 
the index of the $D_6$ curve is $2$
and the index of the critical $D_3$ curve is $3$, the mass formula (\ref{defMG}) gives:
\be
{\chi}_o(\M(\Z_2))=\f{1}{4}+\f{1}{2}-1=-\f{1}{4}
\la{chivZ2}
\ee

\section{Summary and outline}

It is natural to ask  whether the scheme presented here can be improved 
(in particular whether it can be made fully rigorous) or generalized to other functions
on moduli spaces, or to other moduli spaces. 
All these questions are open; however one can make reasonable conjectures on the basis
of this work. First it is natural to expect that any Riemann surface with large automorphism group
is a critical point of any smooth function on the Teichm\"uller space, invariant with respect
to the action of the mapping class group. 
As another confirmation
of this conjecture one can prove this statement for the Klein curve, 
which reads in projective coordinates:
$$
xy^3 +yz^3+zx^3=0\, .
$$
According to \cite{RauchLewittes}, the matrix of $b$-periods of the Klein curve is given by:
$$
{\B}=\begin{pmatrix}
-\frac{1}{8}+\frac{3\sqrt{7}}{8}\,i& -\frac{1}{4}-\frac{\sqrt{7}}{4}\,i & -\frac{3}{8}+\frac{\sqrt{7}}{8}\,i\\
-\frac{1}{4}-\frac{\sqrt{7}}{4}\,i& \frac{1}{2}+\frac{\sqrt{7}}{2}\,i & -\frac{1}{4}-\frac{\sqrt{7}}{4}\,i\\
-\frac{3}{8}+\frac{\sqrt{7}}{8}\,i& -\frac{1}{4}-\frac{\sqrt{7}}{4}\,i & \frac{7}{8}+\frac{\sqrt{3}}{8}\,i
\end{pmatrix}\;;
$$
the matrix
$$
\Sigma=\begin{pmatrix}
1&1&1&1&0&0\\
0&-1&-1&-1&1&0\\
0&1&0&1&-1&0\\
-1&-1&0&0&0&0\\
-1&-1&0&0&0&-1\\
-1&0&0&0&0&-1
\end{pmatrix}
$$
belongs to the stabilizer of ${\B}$ in  $Sp(6,\Z)$.
Then in complete analogy to the proof of Theorem \ref{osnlem} we find 
that the spectrum of the
six-dimensional analog of the matrix ${\mathbb A}$ (\ref{TRANSFORM}) does not contain unity;
thus any smooth  function on the Teichm\"uller  space invariant with respect to the mapping class group has a 
critical point at the 
Klein curve.

Another natural conjecture is that in genera 3 and higher there also exist analogs of the extremal
genus two curve from the $D_3$ family which does not have a large group of automorphisms.
The mass formula for the orbifold Euler characteristic can, perhaps, in some cases 
be helpful in proving the existence of such a curve similarly to the 
genus two case.

A further question is whether the determinant of the Laplacian in the Poincar\'e metric (which is maximal
among all metrics within a given conformal class) has the same set of critical points 
and the same signatures of the Hessian as the determinant in the Bergman metric studied here.
We suppose that the function ${\rm det}\Delta$ 
has only four critical points, similarly to  ${\rm det}\Delta_B$, 
(the existence of three of them we have proved here)
with the maximum at the Burnside curve, but the
fourth critical point from the $D_3$ family is most probably different from the one for 
${\rm det}\Delta_B$. We hope that either 
analytical or numerical analysis of these questions will be possible in the near future.

Finally we hope that the analysis of the global properties of appropriate analogs of ${\rm det}\Delta_B$
should be possible for other interesting spaces:
Hurwitz spaces and the spaces of Abelian and quadratic differentials on Riemann surfaces.
We hope that, say, for spaces of Abelian differentials $w$ on Riemann surfaces \cite{KonZor} 
the proper functional will be the determinant of the Laplacian 
operator in the flat metric with 
conical singularities given by $|w|^2$. Exact formulas for such determinants (and their 
analogs on Hurwitz spaces) obtained in 
\cite{IMRN,DifGeo} should enable at least an efficient numerical analysis of their 
global behavior; that could provide new geometrical information about 
these spaces having been much less  studied than the moduli spaces of Riemann surfaces.

{\bf Acknowledgments}
We are grateful to L.Chekhov, D.Jacobson and P.Zograf for interesting
discussions.
We thank C.Cummins for information about the fundamental domains of
the groups $\Gamma_0(2,3)+$.
The work of DK was partially supported by the Concordia Research Chair 
grant, NSERC and  NATEQ. 
AK thanks the Max Planck institute for Mathematics  in the Sciences in Leipzig
for hospitality and excellent working conditions. DK thanks the Laboratory of Mathematical Physics
of Universit\'e de Bourgogne and the Max-Planck Institute for Mathematics in Bonn, where this work was
completed, for warm hospitality.

\end{document}